\documentclass{amsart}

\usepackage{amsthm}
\usepackage{amsmath}
\usepackage{amssymb}
\usepackage{amsfonts}
\usepackage{dsfont}
\usepackage{mathrsfs}
\usepackage{graphicx}

\usepackage{enumerate}
\usepackage{hyperref}
\usepackage{indentfirst}
\usepackage{color}
\usepackage{multicol,colordvi,epsfig,float}
\usepackage{pgfplots}
\pgfplotsset{width=9cm,compat=1.15}

\newcommand{\norm}[2]{|  {#1}  | _{#2}}
\newcommand{\norma}[2]{\|  {#1} \|_{#2}}
\newcommand{\prodl}[2]{ (\ {#1} \ , \ {#2} \ )}

\newcommand{\D}[1]{\displaystyle{#1}}

\renewcommand{\qedsymbol}{$\blacksquare$}

\usepackage{verbatim}

\newtheorem{theorem}{Theorem}[section]
\newtheorem{lemma}[theorem]{Lemma}
\newtheorem{corollary}[theorem]{Corollary}

\theoremstyle{definition}
\newtheorem{definition}[theorem]{Definition}

\theoremstyle{remark}
\newtheorem{remark}[theorem]{\bf Remark}

\numberwithin{equation}{section}

\begin{document}

\title[Optimal control for a chemotaxis-consumption model]{
Optimal control related to weak solutions of a chemotaxis-consumption model}

\author{Francisco Guillén-González}
\address{Departamento de Ecuaciones Diferenciales y Análisis Numérico, Universidad de Sevilla}
\email{guillen@us.es}

\author{André Luiz Corrêa Vianna Filho}
\address{Departamento de Ecuaciones Diferenciales y Análisis Numérico, Universidad de Sevilla}
\email{acorreaviannafilho@us.es}





\begin{abstract}
  In the present work we investigate an optimal control problem related to the following chemotaxis-consumption model in  a bounded domain $\Omega\subset \mathbb{R}^3$:
  $$\partial_t u - \Delta u  = - \nabla \cdot (u \nabla v), \quad \partial_t v - \Delta v  = - u^s v + f \,v\, 1_{\Omega_c},$$
  with $s \geq 1$, endowed with isolated boundary conditions and initial conditions for $(u,v)$, being $u$ the cell density, $v$ the chemical concentration and $f$ the control acting in the $v$-equation through the bilinear term $f \,v\, 1_{\Omega_c}$, in a subdomain $\Omega_c \subset \Omega$. We address the existence of optimal control restricted to a weak solution setting, where, in particular, uniqueness of state $(u,v)$ given a control $f$ is not clear. Then by considering weak solutions satisfying an adequate  energy inequality, we prove the existence of optimal control subject to uniformly bounded controls. Finally, we discuss the relation between the considered control problem and two other related ones, where the existence of optimal solution can not be proved.
\end{abstract}

\maketitle



\section{Introduction}

  We consider an optimal control problem subject to a chemotaxis-consumption model. Let $\Omega\subset \mathbb{R}^3$ be a bounded domain, denoting by $\Gamma$ its boundary, and define $Q: = (0,T) \times \Omega$, for a given final time $T > 0$. Let $u=u(t,x)$ and $v=v(t,x)$ be the density of cell population and the concentration of chemical substance, respectively, defined in $(t,x) \in Q$. The aforementioned chemotaxis-consumption model is given by the PDE system
  \begin{equation} \label{problema_P}
    \left\{\begin{array}{l}
      \partial_t u - \Delta u  = - \nabla \cdot (u \nabla v), \quad
      \partial_t v - \Delta v  = - u^s v, \quad \hbox{in $Q$,}\\ [6pt]
      \partial_\eta u  |_{\Gamma}  =  \partial_\eta v  |_{\Gamma} = 0, \quad
      u(0,\cdot)  = u^0, \quad v(0, \cdot) = v^0, \quad \hbox{in $\Omega$,}
    \end{array}\right.
  \end{equation}
  where $s \geq 1$, $\partial_\eta u  |_{\Gamma}$ denotes the outer normal derivative of $u$ on the boundary and  $u^0, v^0 \geq 0$ in $\Omega$ are the initial conditions. 
  
 The controlled problem consists of the chemotaxis-consumption model \eqref{problema_P} with a term $f v 1_{\Omega_c}$ added to the chemical equation:
  \begin{equation} \label{problema_P_controlado}
    \left\{\begin{array}{l}
      \partial_t u - \Delta u  = - \nabla \cdot (u \nabla v), \quad
      \partial_t v - \Delta v  = - u^s v + f v 1_{\Omega_c},
      \quad \hbox{in $Q$,} \\ [6pt]
      \partial_\eta u  |_{\Gamma}  =  \partial_\eta v  |_{\Gamma} = 0, \quad
      u(0)  = u^0, \quad v(0) = v^0, \quad \hbox{in $\Omega$,}
    \end{array}\right.
  \end{equation}
  where $\Omega_c \subset \Omega$ is the control subdomain and $f: Q_c:=(0,T)\times \Omega_c \rightarrow \mathbb{R}$ is the control function. This term $fv1_{\Omega_c}$ means that we are controlling the system by acting directly on the chemical equation on the subdomain $\Omega_c$, increasing (where $f\ge 0$) or decreasing  (where $f\le 0$) the chemical concentration, while the control on the cell equation is exerted indirectly. The use of this   bilinear term $fv1_{\Omega_c}$ as control allows to preserve the positivity of $v$ independently of the sign of $f$. The control of chemotaxis systems through the direct action on the chemical equation has been considered in previous studies \cite{ryu2001optimal,Fister-Mccarthy,guillen2020optimal,guillen2020bi,guillen2020regularity,lopez2021optimal,silva2022bilinear,tang2022optimal,guillen2023optimal}. 
  
  To contextualize this work, we highlight the closest related previous literature, beginning with some theoretical studies on the uncontrolled problem \eqref{problema_P}. For the case $s = 1$,  the existence of a global weak solution in smooth and convex 3D domains is proved in \cite{tao2012eventual}. Moreover, these solutions become smooth after a large enough period of time and their long time behavior is studied. In \cite{tao2019global}, a parabolic-elliptic version of \eqref{problema_P} is studied, leading to results on the existence, uniqueness and long time behavior of a global classical solution in $n$-dimensional smooth domains.
  
  Still considering $s = 1$, some authors have also focused on the model \eqref{problema_P} but coupled with a surrounding fluid (called chemotaxis-Navier-Stokes model), as a particular case of the models introduced by \cite{tuval2005bacterial}. In \cite{winkler2012global}, the existence of global weak solutions is proved in smooth and convex $3D$ domains. Moreover, these solutions become smooth after a large enough period of time and their long time behavior is studied. These results are extended in \cite{jiang2015global}  to nonconvex but smooth domains. In \cite{winkler2014stabilization} the author studies the asymptotic behavior of the chemotaxis-Navier-Stokes equations in 2D domains for adequate generalizations of the chemotaxis and consumption terms, proving the convergence towards constant states in the $L^{\infty}$-norm. In \cite{winkler2016global}, existence of global weak solutions for the chemotaxis-Navier-Stokes equations is established in $3$D smooth and convex domains, and the study of the asymptotic behavior of these solutions is carried out in \cite{winkler2017far}.

  All the aforementioned works are related to the chemotaxis-consumption model \eqref{problema_P} with $s = 1$, i.e. the chemotaxis and consumption terms are given, respectively, by
  \begin{equation} \label{chemotaxis_consumption_terms_s=1}
    \nabla \cdot (u \nabla v) \mbox{ and } -uv.
  \end{equation}
  However, from the modeling point of view, it may be useful to study models with more general terms. Indeed, different assumptions during the modeling process can lead to different terms in the equations and, therefore, it would be useful to have results for a model with more general terms, rather than for a specific case. In particular, when the potential consumption $-u^s v$ is considered, the action of the consumption term is more accentuated because, compared to the bilinear case $-u v$, the consumption is stronger when $u > 1$ and is weaker for $u < 1$. Moreover, it is also interesting to know, for instance, what is the effect of different chemotaxis and consumption terms in the properties of the solutions, such as regularity, boundedness, asymptotic behavior and so on.

  In \cite{winkler2012global}, already cited above, the chemotaxis and consumption terms considered had the form $\nabla \cdot (u F'(u) \chi(v) \nabla v)$ and $-F(u)h(v)$, including \eqref{chemotaxis_consumption_terms_s=1} as a particular case. The results were then obtained under hypotheses on the behavior of the functions $\chi(v)$, $h(v)$, $F(u)$ and its derivative $F'(u)$. More recently, problem \eqref{problema_P} has been studied in \cite{ViannaGuillen2023uniform} and the effects of any consumption power $s\ge 1$ in the regularity of the solutions have been addressed. In addition, the class of considered domains has been enlarged, including $2$D and $3$D domains which are not necessarily smooth or convex.
  
  Another relevant topic in chemotaxis models is the existence or non-existence of blowing-up solutions. If we consider model \eqref{problema_P}, this question has been answered in $2D$ domains. In fact, considering $s = 1$ and $2D$ convex domains, existence and uniqueness of a global classical solution that is uniformly bounded up to infinity time is proved in \cite{tao2012eventual}; and for $s \geq 1$ and more general nonconvex $2D$ domains, existence and uniqueness of a global strong solution that is uniformly bounded up to infinity time is proved in \cite{ViannaGuillen2023uniform}. As far as we know, this question remains open for $3D$ domains. In fact, some authors investigated hypotheses that could lead to no blow-up results in $n$-dimensional spaces, for $n \geq 3$. They obtained advances under smallness assumptions on the chemotaxis coefficient and on  $\norma{v^0}{L^{\infty}(\Omega)}$. On this subject, we refer the interested reader to \cite{tao2011boundedness} and \cite{baghaei2017boundedness}, for the problem \eqref{problema_P} with $s = 1$. In addition,  these results are extended to other related chemotaxis models with consumption in  \cite{fuest2019analysis} and \cite{frassu2021boundedness}.
  
  The works cited so far address the analysis of the chemotaxis-consumption model \eqref{problema_P} and helps us understand how the system evolves from the given initial data. However, when it comes to PDE models which describe physical phenomena, such as the chemotaxis models, as important as the analysis of the system itself, are the studies of control problems related to them. Concerning chemotaxis models, particularly relevant are the optimal control problems. Due to the low regularity of the weak solutions and the lack of uniqueness results in $3D$, some works are concentrated in $2D$ domains, where one usually has existence and uniqueness of strong solution to the controlled problem, which allows to prove the existence of optimal control, and to derive an optimality system, where existence and regularity of Lagrange multipliers are deduced. In this direction, we refer the reader to the works on control problems related to: the Keller-Segel model \cite{ryu2001optimal}; a chemorepulsion-production model \cite{guillen2020optimal,guillen2020bi}; a Keller-Segel logistic model \cite{silva2022bilinear}; a chemotaxis model with indirect consumption \cite{yuan2022optimal}; and a chemotaxis-haptotaxis model \cite{tang2022optimal}.
  
  For optimal control problems related to chemotaxis models in $3D$ domains the situation is more complex because, in general, we have results of existence of weak solutions, however, in many cases, there is not any result on the existence and uniqueness  of global in time strong solutions. In this setting, it is usual to introduce a regularity criterion, which is a mild additional regularity hypothesis, sufficient to conclude that a weak solution satisfying this regularity is actually the unique strong solution. A didactic introduction to this kind of adaptation can be found in \cite{casas1998optimal}, for  an optimal control problem related to the Navier-Stokes equations in $3D$ domains. We refer the reader to \cite{guillen2020regularity}, where a regularity criterion is established  to study an optimal  control problem related to a chemorepulsion-production model in $3D$ domains. The drawback of using a regularity criterion is that it is not clear in general if the admissible set  is nonempty. In \cite{guillen2020regularity}, the authors show that if $\Omega_c = \Omega$, that is, if the control acts in the whole domain, and the initial chemical density $v^0$ is strictly positive and separated from zero, then the admissible set is nonempty. To do that, the idea is to define the control $f$ \emph{a posteriori}, depending on a regular pair $(u,v)$.
  
\

  Although the study of optimal control problems related to chemotaxis models is an interesting and growing topic, we still have a relative low number of studies concerning the chemotaxis-consumption model \eqref{problema_P}. In fact, to the best of our knowledge, we can cite two works: \cite{lopez2021optimal} and \cite{guillen2023optimal}. In \cite{lopez2021optimal}, the authors apply a regularity criterion to a chemotaxis-consumption-Navier-Stokes model, proving the existence of optimal control and the existence of Lagrange multipliers, arriving at the optimality system. In \cite{guillen2023optimal}, a controlled problem subject to strong solutions of \eqref{problema_P_controlado} is studied. A sharp regularity criterion is proved and used to show the existence of optimal solution and the existence and uniqueness of Lagrange multipliers associated to any local optimal solution. A generic PDE linear system is also used to study the regularity of the Lagrange multipliers. In both \cite{lopez2021optimal} and \cite{guillen2023optimal}, it is only possible to prove that the admissible set is nonempty for controls acting in the whole domain $\Omega$.
  
  In view of all exposed so far, the objective of the present work is to study an optimal control problem related to \eqref{problema_P} for which we are able to prove the existence of global optimal solution in the weak solutions setting, that is, without using any regularity criterion or hypothesis over the admissible set. To achieve it, an important novelty of this work is to consider weak solutions of the controlled model \eqref{problema_P_controlado} satisfying an energy inequality (see \eqref{desigualdade_integral_limite} below). Afterwards, we prove existence of global optimal solution and, to conclude, we discuss the relation between this optimal control problem and two other related ones, where the existence of optimal solution can not be proved.
  
  We remark that the deduction of first order optimality conditions in this weak solution setting remains as an open problem. Indeed, since we only have the weak regularity, it is not possible to prove the well-posedness of the linearized problem around a local optimal solution, which is the essential hypothesis to apply a generic Lagrange multiplier method as in \cite{guillen2023optimal}. Also, given a control $f$, we do not have in general uniqueness of the state $(u,v)$, hence it is not possible to define the  ``control-to-state'' mapping and follow the procedure used in \cite{lopez2021optimal} to compute the derivative of the state with respect to the control.
  

  \subsection{Main results}
  
  Along this paper, we impose the following hypotheses:
   \begin{equation*}
   \hbox{$\Omega \subset \mathbb{R}^3$ is a bounded domain with boundary $\Gamma$ of class $C^{2,1}$,}
     \end{equation*}
  \begin{equation} \label{condicao_sobre_o_controle_f}
    f \in L^q(Q_c), \mbox{ for some } q > 5/2,
  \end{equation}
  \begin{equation} \label{cond_ini_propriedades}
    (u^0, v^0) \in L^p(\Omega) \times W^{2-2/q,q}(\Omega),
  \end{equation}
  with $p = 1 + \varepsilon$, for some $\varepsilon > 0$, if $s = 1$, and $p = s$, if $s > 1$.
  
  Let us define the specific functional spaces appearing for the weak solution setting. For $s \in [1,2)$,
  \begin{equation*}
    \begin{array}{rl}
      X_u = & \hspace{-4pt} \Big \{ u \ | \ u\in L^{\infty}(0,T;L^s(\Omega)) \cap  L^{5s/3}(Q), \  \nabla u \in   L^{5s/(3 + s)}(Q), \\
      & \qquad \partial_t u \in L^{5s/(3 + s)}(0,T;W^{1,5s/(4s-3)}(\Omega)') \Big \},
    \end{array}
  \end{equation*}
  for $s \geq 2$,
  \begin{equation*}
    \begin{array}{rl}
      X_u = & \Big \{ u\ | \ u \in L^{\infty}(0,T;L^s(\Omega)) \cap L^{5s/3}(Q),\  \nabla u \in L^{2}(Q), \\
      & \qquad  \partial_t u \in L^2(0,T;H^1(\Omega)') \Big \},
    \end{array}
  \end{equation*}
 and  for $s \geq 1$,
  \begin{equation*}
    \begin{array}{rl}
      X_v = & \hspace{-4pt} \Big \{ v\ | \ v \in L^{\infty}(Q), \ \nabla v \in L^{\infty}(0,T;L^2(\Omega))  \cap L^4(Q) \cap  L^2(0,T;H^1(\Omega)), \\
      & \qquad \Delta v \in L^2(Q), \ \partial_t v \in L^{5/3}(Q) \Big \}.
    \end{array}
  \end{equation*}
  We also introduce  the bounded convex set for the control
  \begin{equation*}
    B_q(M) = \Big \{ f \in L^q(Q_c) \ | \ \norma{f}{L^q(Q_c)} \leq M \Big \}.
  \end{equation*}
  
  \begin{definition}[\bf Weak Solution of \eqref{problema_P_controlado}]
    A pair $(u,v)$ is called a weak solution of \eqref{problema_P_controlado} if $u(t,x),v(t,x) \geq 0$ $a.e.$ $(t,x) \in Q$, with
    \begin{equation*}
      u \in X_u, \ v \in X_v
    \end{equation*}
    and satisfying the initial conditions for $(u,v)$, the $u$-equation of \eqref{problema_P_controlado} and the boundary condition of $u$ in the variational sense
    \begin{equation*}
      \langle\partial_t u , \varphi \rangle + \int_\Omega \nabla u\cdot\nabla \varphi \ dx  = \int_\Omega u \nabla v\cdot \nabla \varphi \ dx,
      \quad a.e. \ t > 0, 
      \quad \forall \, \varphi \in 
      W^{1,5s/(4s-3)}(\Omega),
    \end{equation*}
     the $v$-equation of \eqref{problema_P_controlado} point-wisely $a.e.$ $(t,x) \in Q$ (in fact, the $v$-equation is satisfied in $L^{5/3}(Q)$) and, since $\Delta v \in L^2(Q)$, the boundary condition of $v$ in the sense of $H^{-1/2}(\Gamma)$.
    \hfill $\square$
  \end{definition}

  \begin{remark} \label{remark_initial_conditions_weak_solutions}
    Considering the regularity of  the weak solution pair $(u,v)$ and the regularity of the time derivatives, $u_t$ and $v_t$, we can specify the sense in which the initial conditions are attained. Indeed, 
    we are able to conclude that $(u,v)$ is weakly continuous from $[0,\infty)$ into $L^s(\Omega) \times H^1(\Omega)$, if $s \in [1,2]$, and $L^2(\Omega) \times H^1(\Omega)$, if $s \geq 2$ (see \cite[Remark 2]{ViannaGuillen2023uniform}).
    \hfill $\square$
  \end{remark}

  \
  
  The proof of existence of weak solution to the controlled problem \eqref{problema_P_controlado} is based in the treatment of the uncontrolled model \eqref{problema_P} given in \cite{ViannaGuillen2023uniform} and extended to the model \eqref{problema_P_controlado}, which has a non-smooth control $f$ as a coefficient. An important step in \cite{ViannaGuillen2023uniform} is the obtaining of an energy inequality using the change of variable from $(u,v)$ to $(u,z)$, with $z = \sqrt{v + \alpha^2}$, where $\alpha > 0$ is a sufficiently small but fixed real number, independently of $(u,v)$ and $f$,  which will be  chosen in Lemmas \ref{lemma_estimativa_u_v_m_1_s=1}, \ref{lemma_estimativa_u_v_m_1_s_intermediario} and \ref{lemma_estimativa_u_v_m_1_s_geq_2} below. Here, the energy inequality  satisfied by the constructed weak solutions of \eqref{problema_P_controlado} will also be written in terms of $(u,z)$. In fact, we consider the energy
  \begin{equation*}
    E(u,z)(t) = \frac{s}{4} \D{\int_{\Omega}}{g(u(t,x)) \ dx} + \frac{1}{2} \int_{\Omega}{\norm{\nabla z(t,x)}{}^2} \ dx,
  \end{equation*}
  where
  \begin{equation*}
    g(u) = \left \{ 
    \begin{array}{rl}
      (u+1)ln(u+1) - u, & \mbox{if } s = 1,  \\[3pt]
      \dfrac{1}{s(s-1)} \, u^s, & \mbox{if } s > 1. 
    \end{array}
    \right.
  \end{equation*}
    
\
  
  We have the following result of existence of weak solutions to \eqref{problema_P_controlado}.
  
  \begin{theorem}[\bf Existence of energy inequality weak solutions] \label{teo_existencia_solucao_fraca}
    Given $f \in L^q(Q)$ ($q > 5/2$), there is a non-negative weak solution $(u,v)$ of \eqref{problema_P_controlado} satisfying the following energy inequality
    \begin{equation} \label{desigualdade_integral_limite}
      \begin{array}{l}
        E(u,z)(t_2) + \beta \D{\int_{t_1}^{t_2} \int_{\Omega}}{\norm{\nabla [u + 1]^{s/2}}{}^2 \ dx} \ dt + \dfrac{1}{4} \D{\int_{t_1}^{t_2} \int_{\Omega}}{u^s \norm{\nabla z}{}^2 \ dx} \ dt \\[6pt]
        + \beta \Big ( \D{\int_{t_1}^{t_2} \int_{\Omega}}{\norm{D^2 z}{}^2 \ dx} \ dt + \D{\int_{t_1}^{t_2} \int_{\Omega}}{\frac{\norm{\nabla z}{}^4}{z^2} \ dx} \ dt \Big ) \\[10pt]
        \leq E(u,z)(t_1) + \mathcal{K}(\norma{f}{L^q(Q)}, \norma{v_0}{W^{2-2/q,q}(\Omega)}),
      \end{array}
    \end{equation}
    for $a.e.$ $t_1,t_2 \in [0,T]$. Here, $\mathcal{K}(\norma{f}{L^q(Q)},\norma{v_0}{W^{2-2/q,q}(\Omega)})$ is a continuous and increasing function with respect to $\norma{f}{L^q(Q)}$ and $\beta > 0$ is a constant, independent of $(u,v,f)$. Moreover, inequality \eqref{desigualdade_integral_limite} is also valid for $t_1 = 0$, with $E(u,z)(0) = E(u^0,v^0)$. Finally, in the case $s > 1$, inequality \eqref{desigualdade_integral_limite} also holds  for all $t_2 \in (t_1,T]$.
  \end{theorem}
  
  \begin{remark}
    The existence of weak solutions satisfying an energy inequality is commonly seen, for instance, for fluid models, and is used to prove either weak-strong uniqueness results \cite{lions1996incompressible} or large  time behaviour  \cite{miyakawa1988energy}. In the present work, we use this ``energy inequality weak solution" setting in order to prove the existence of global optimal solution. To the best of our knowledge, this is the first time that the concept of weak solution with energy inequality  is applied to this purpose for chemotaxis models.
    \hfill $\square$
  \end{remark}
  
  \
  
  Next we introduce the minimization problem. Consider the functional 
  $$J: L^{5s/3}(Q) \times L^2(Q) \times L^q(Q) \longrightarrow \mathbb{R}$$ given by
  \begin{align*}
    & J(u,v,f) : = \dfrac{3\gamma_u}{5s} \int_0^T{\norma{u(t) - u_d(t)}{L^{5s/3}(\Omega)}^{5s/3} \ dt} \\
    & \D + \dfrac{\gamma_v}{2} \int_0^T{\norma{v(t) - v_d(t)}{L^2(\Omega)}^2 \ dt} + \dfrac{\gamma_f}{q} \int_0^T{\norma{f(t)}{L^q(\Omega)}^q \ dt},
  \end{align*}
  where $(u_d,v_d) \in L^{5s/3}(Q) \times L^2(Q)$ represents the desired states, $\gamma_u, \gamma_v > 0$ measure the errors in the states and $ \gamma_f>0$ the cost of the control. In view of the existence result, Theorem \ref{teo_existencia_solucao_fraca}, one could expect the following admissible sets
  \begin{equation*}
    \begin{array}{rl}
      S_{ad}^w = & \hspace{-4pt} \Big \{ (u,v,f) \in X_u \times X_v \times L^q(Q) \ | \ (u,v) \mbox{ is a} \\
      & \quad \mbox{weak solution of \eqref{problema_P_controlado} with control } f \}
    \end{array}
  \end{equation*}
  or
  \begin{equation*}
    \begin{array}{rl}
      S_{ad}^E = & \hspace{-4pt} \Big \{ (u,v,f) \in X_u \times X_v \times L^q(Q) \ | \ (u,v) \mbox{ is a weak solution of \eqref{problema_P_controlado}} \\
      & \quad \mbox{with control } f \mbox{ and satisfies the energy inequality } \eqref{desigualdade_integral_limite} \Big \}
    \end{array}
  \end{equation*}
  and then state the corresponding minimization problems
  \begin{equation} \label{problema_de_minimizacao_S_ad^w}
    \left \{
    \begin{array}{l} \D
      min \ J(u,v,f) \\
      \mbox{subject to } (u,v,f) \in S_{ad}^w,
    \end{array}
    \right.
  \end{equation}
  or
  \begin{equation} \label{problema_de_minimizacao_S_ad^E}
    \left \{
    \begin{array}{l} \D
      min \ J(u,v,f) \\
      \mbox{subject to } (u,v,f) \in S_{ad}^E.
    \end{array}
    \right.
  \end{equation}
  Thanks to Theorem \ref{teo_existencia_solucao_fraca} we have that both $S_{ad}^w$ and $S_{ad}^E$ are nonempty sets. However, we are not able to prove that problem \eqref{problema_de_minimizacao_S_ad^w} or \eqref{problema_de_minimizacao_S_ad^E} has a solution, as we will analyze in Remark \ref{remark_open_problems} and Subsection \ref{subsec: existencia_controle_otimo}, respectively.
  
  Therefore, in order to find an optimal control related to weak solutions of \eqref{problema_P_controlado}, we introduce the following admissible set, for each $M>0$:
  \begin{equation*}
    \begin{array}{rl}
      S_{ad}^M = & \hspace{-4pt} \Big \{ (u,v,f) \in X_u \times X_v \times B_q(M) \ | \ (u,v) \mbox{ is a weak solution of}  \\
      & \mbox{\eqref{problema_P_controlado} with control } f \mbox{ and satisfies \eqref{desigualdade_integral_limite} changing the} \\
      & \mbox{constant } \mathcal{K}(\norma{f}{L^q(Q)}, \norma{v_0}{W^{2-2/q,q}(\Omega)}) \mbox{ by } \mathcal{K}(M, \norma{v_0}{W^{2-2/q,q}(\Omega)}) \Big \}
    \end{array}
  \end{equation*}
  and the corresponding minimization problem
  \begin{equation} \label{problema_de_minimizacao}
    \left \{
    \begin{array}{l} \D
      min \ J(u,v,f) \\
      \mbox{subject to } (u,v,f) \in S_{ad}^M.
    \end{array}
    \right.
  \end{equation}
  Again, from Theorem \ref{teo_existencia_solucao_fraca}, one has $S_{ad}^M \neq \emptyset$. But now, we are able to prove the following result.
  \begin{theorem}[\bf Existence of optimal control] \label{teo_existencia_controle_otimo}
    For each $M > 0$, the optimal control problem \eqref{problema_de_minimizacao} has at least one global optimal solution, that is, there is $(\overline{u},\overline{v},\overline{f}) \in S_{ad}^M$ such that
    \begin{equation*}
      J(\overline{u},\overline{v},\overline{f}) = \min_{(u,v,f) \in S_{ad}^M} J(u,v,f).
    \end{equation*}
  \end{theorem}

\

  \begin{remark} \label{remark_energy_structure}
    As it can be observed in Subsection \ref{subsec: existencia_controle_otimo}, to prove Theorem \ref{teo_existencia_controle_otimo}, it will be fundamental to have an energy structure such as inequality \eqref{desigualdade_integral_limite}. In fact, considering the minimizing sequence argument used in Subsection \ref{subsec: existencia_controle_otimo}, the energy inequality \eqref{desigualdade_integral_limite} is the key point to guarantee that all the possible limits of the minimizing sequence are weak solutions of the controlled model. In fact, the corresponding energy estimates must be strong enough to guarantee that the possible limits of the minimizing sequence are weak solutions of the controlled model. Therefore, if the model does not admit an energy structure, as it seems to be the case in \cite{winkler2015large}, for example, it is not clear how to prove the existence of optimal solution.
    \hfill $\square$
  \end{remark}

\
  
  By construction, we have the following relation between problems \eqref{problema_de_minimizacao_S_ad^w} and \eqref{problema_de_minimizacao}:
  \begin{equation}\label{w-M}
    J_{inf}^w:=\inf_{(u,v,f) \in S_{ad}^w} J(u,v,f) \leq \min_{(u,v,f) \in S_{ad}^M} J(u,v,f).
  \end{equation}
  On the other hand, in this paper we will obtain the following relation between the minimization problems \eqref{problema_de_minimizacao_S_ad^E} and \eqref{problema_de_minimizacao} for $M$ large enough:
  \begin{theorem} \label{teo_relacao_entre_problemas_de_minimizacao}
    If
    \begin{equation*}
      M^q \geq \frac{q}{\gamma_f} \inf_{(u,v,f) \in S_{ad}^E}{J(u,v,f)},
    \end{equation*}
    we have the inequality
    \begin{equation}\label{M-E}
      \min_{(u,v,f) \in S_{ad}^M} J(u,v,f) \leq \inf_{(u,v,f) \in S_{ad}^E}{J(u,v,f)}:=J_{inf}^E.
    \end{equation}
  \end{theorem}
  
  \
   
  \begin{remark} \label{remark_open_problems}
    From Theorem \ref{teo_existencia_controle_otimo}, for each $M > 0$ there is $(u^M,v^M,f^M) \in S_{ad}^M$ such that
    \begin{equation*}
      J(u^M,v^M,f^M) = \min_{(u,v,f) \in S_{ad}^M} J(u,v,f).
    \end{equation*}
    Let $M_2 > M_1 > 0$. Since $S_{ad}^{M_1} \subset S_{ad}^{M_2}$ then  $J(u^M,v^M,f^M)$ decreases as $M$ increases. Therefore, since $J(u^M,v^M,f^M)$ is bounded from below, there exists $\D{\lim_{M \to \infty} \ J(u^M,v^M,f^M)}$ and, accounting for \eqref{w-M} and \eqref{M-E}, one has the inequalities 
    \begin{equation*}
      J_{inf}^w   \leq \lim_{M \to \infty} \ J(u^M,v^M,f^M) \leq  J_{inf}^E.
    \end{equation*}
    Let $(u^{\infty},v^{\infty},f^{\infty}) \in L^{5s/3}(Q) \times L^2(Q) \times L^q(Q)$ be the weak limit of a subsequence of $\{ (u^M,v^M,f^M) \}_M$. Then, the weakly lower semicontinuity of $J$ in $L^{5s/3}(Q) \times L^2(Q) \times L^q(Q)$  leads  to
    \begin{equation} \label{J^infty}
      J(u^{\infty},v^{\infty},f^{\infty}) \leq \lim_{M \to \infty} \ J(u^M,v^M,f^M).
    \end{equation}
    
    In our opinion, the proof or the refutation of the following two questions are interesting open problems:
    \begin{enumerate}
      \item  $(u^{\infty},v^{\infty},f^{\infty}) \in S_{ad}^w$ ? \label{open_prob_1}
      \item  $\D{\lim_{M \to \infty}} \ J(u^M,v^M,f^M) = J_{inf}^w$ ? \label{open_prob_2} 
    \end{enumerate}
    In fact, if (\ref{open_prob_2}) were valid, then  $J_{inf}^w$ could be approximated by $\D{\min_{(u,v,f) \in S_{ad}^M}} J(u,v,f)$ as $M \to \infty$. On the other hand, if (\ref{open_prob_1}) and (\ref{open_prob_2}) were valid, then $(u^{\infty},v^{\infty},f^{\infty})$ becomes an optimal solution of \eqref{problema_de_minimizacao_S_ad^w}. Indeed, from (\ref{open_prob_1}) we have $J_{inf}^w \leq J(u^{\infty},v^{\infty},f^{\infty})$, and from (\ref{open_prob_2}) and \eqref{J^infty}, we have $J(u^{\infty},v^{\infty},f^{\infty}) \leq J_{inf}^w$.
    \hfill $\square$
  \end{remark}
  
  The rest of the paper is organized as follows. In Section \ref{Se:Preliminary} we present some preliminary results. The existence of weak solutions satisfying the energy inequality \eqref{desigualdade_integral_limite} for the controlled problem (proof of Theorem \ref{teo_existencia_solucao_fraca}) is established in Section \ref{Se:Existence} and, in Section \ref{section:existence of optimal control} we study the optimal control problem, proving Theorems \ref{teo_existencia_controle_otimo} and \ref{teo_relacao_entre_problemas_de_minimizacao}.
  

\section{Preliminary results}\label{Se:Preliminary}
  
  \begin{lemma} \label{interpolacao_norma_10/3}
    Let $\Omega \subset \mathbb{R}^3$ be a bounded Lipschitz domain. We have
    \begin{equation*}
      \norma{v}{L^{10/3}(\Omega)} \leq C \norma{v}{L^2(\Omega)}^{2/5} \norma{v}{H^1(\Omega)}^{3/5}.
    \end{equation*}
  \end{lemma}
  
  \begin{lemma}{\bf (\cite{brezis2011functional})} \label{lemma_semicontinuidade_fraca_inferior}
    Let B be a Banach space, let $\{ w_n \}$ be a sequence in $B$ and $w \in B$. Either if $w_n \rightarrow w$ weakly* or weakly in $B$ then $\{ w_n \}$ is bounded in $B$ and $\norma{w}{B} \leq \liminf \norma{w_n}{B}$.
  \end{lemma}
  
  \begin{lemma}{\bf (\cite{Temam})} \label{lemma_continuidade_fraca}
    Let $X$ and $Y$ be two Banach spaces such that $X \subset Y$ with a continuous injection. If $\phi \in L^{\infty}(0,T;X)$ and $\phi \in C([0,T];Y)$, then $\phi \in C_w([0,T];X)$.
  \end{lemma}
  
  \begin{lemma}{\bf (\cite{feireisl2009singular})} \label{lema_regularidade_eq_calor}
    Let $\Omega$ be a bounded domain of $\mathbb{R}^N$ such that $\Gamma$ is of class $C^2$. Let $p \in (1,3)$, $w^0 \in W^{2-2/p,p}(\Omega)$ and $h \in L^p(Q)$. Then the problem
    \begin{equation*}
      \left\{\begin{array}{l}
        \partial_t w - \Delta w = h, \quad\mbox{in } Q, \\
        \partial_\eta w |_{\Gamma}  =  0, \quad\mbox{on } (0,T) \times \Gamma, \\
        w(0,\cdot)  = w^0, \quad\mbox{in } \Omega,
      \end{array}\right.
    \end{equation*}
    has a unique solution 
    \begin{equation*}
      w \in C([0,T]; W^{2-2/p,p}(\Omega)) \cap L^p(0,T;W^{2,p}(\Omega)), \ \partial_t w \in L^p(Q).
    \end{equation*}
    Moreover, there is a positive constant $C = C(p,T,\Omega)$ such that
    \begin{equation} \label{dependencia_continua_de_w_em_h}
      \begin{array}{l}
        \norma{w}{C([0,T];W^{2 - 2/p,p}(\Omega))} + \norma{w}{L^p(0,T;W^{2,p}(\Omega))} + \norma{\partial_t w}{L^p(Q)} \\[8pt]
        \leq C ( \norma{h}{L^p(Q)} + \norma{w^0}{W^{2 - 2/p,p}(\Omega)}).
      \end{array}
    \end{equation}
  \end{lemma}
  
  \begin{lemma}[\cite{ViannaGuillen2023uniform}] \label{lema_T^m_elevado_a_s}
    Let $w_1$ and $w_2$ be nonnegative real numbers. For each $s \geq 1$ we have
    \begin{equation*}
      \norm{w_2^s - w_1^s}{} \leq s \norm{w_2 + w_1}{}^{s-1} \norm{w_2 - w_1}{}.
    \end{equation*}
  \end{lemma}
    
  Using Lemma \ref{lema_T^m_elevado_a_s}, we can prove the following.
  
  \begin{lemma} \label{lema_convergencia_w_elevado_a_s}
    Let $p \in (1, \infty)$ and let $\{ w_n \}$ be a sequence of nonnegative functions in $L^p(Q)$ such that $w_n \to w$ in $L^p(Q)$ as $n \to \infty$. Then, for every $r \in (1,p)$, $w_n^r \to w^r$ in $L^{p/r}(Q)$ as $n \to \infty$.
  \end{lemma}
  
  \begin{lemma}{\bf (Compactness in Bochner spaces \cite{Simon1986compact})} \label{lema_Simon}
    Let $X,B$ and $Y$ be Banach spaces, let
    \begin{equation*}
      F \subset \Big \{ f \in L^1(0,T;Y) \ | \ \partial_t f \in L^1(0,T;Y) \Big \}.
    \end{equation*}
    Suppose that $X \subset B \subset Y$, with compact embedding $X \subset B$ and continuous embedding $B \subset Y$. Let the set $F$ be bounded in $L^q(0,T;B) \cap L^1(0,T;X)$, for $1 < q \leq \infty$, and $\Big \{ \partial_t f, \ \forall f \in F \Big \}$ be bounded in $L^1(0,T;Y)$. Then $F$ is relatively compact in $L^p(0,T;B)$, for $1 \leq p < q$.
  \end{lemma}
  
  
\


\section{Existence of the controlled problem}\label{Se:Existence}
  
  The existence of weak solutions of the uncontrolled problem \eqref{problema_P} is proved in \cite{ViannaGuillen2023uniform}, by means of a sequence of truncated problems. To define these truncated problems, we are going to use a mollifier regularization of the control $f \in L^q(Q)$ defined via convolution,  considering a sequence  
  (see \cite{brezis2011functional})
  \begin{equation} \label{properties_f_m}
    \begin{array}{c}
      f_m \in C^{\infty}_c(\overline{Q}), \quad \norma{f_m}{L^r(Q)} \leq \norma{f}{L^r(Q)}, \, \mbox{for } r \in [1,q], \\[6pt]
      f_m \rightarrow f \mbox{ strongly in } L^q(Q).
    \end{array}
  \end{equation}
  
  Then, we prove the existence of solution of the controlled problem \eqref{problema_P_controlado} satisfying, in addition, the energy inequality \eqref{desigualdade_integral_limite}, using the following controlled truncated problems:
  \begin{equation}
    \left\{
    \begin{array}{l}
      \partial_t u_m - \Delta u_m  = - \nabla \cdot (T^m(u_m) \nabla v_m), \\
      \partial_t v_m - \Delta v_m  = - T^m(u_m)^s v_m + f_m v_m 1_{\Omega_c}, \\
      \partial_\eta u_m |_{\Gamma}  = \partial_\eta v_m |_{\Gamma} = 0, \quad
      u_m(0)  = u^0_m, \quad v_m(0) = v^0,
    \end{array}
    \right.
    \label{problema_P_m_controlado}
  \end{equation}
  for each $m \in \mathbb{N}$, where the truncation function $T^m \in C^2(\mathbb{R})$ is defined by
  \begin{equation*} 
    T^m(u) = \left \{
    \begin{array}{cl}
      - 1, & \mbox{ if } u \leq -1,  \\
      C^2 \mbox{ extension}, & \mbox{ if } u \in (-1,0), \\
      u, & \mbox{ if } u \in [0, m], \\
      C^2 \mbox{ extension}, & \mbox{ if } u \in (m,m + 2), \\
      m + 1, & \mbox{ if } u \geq m + 2,
    \end{array}
    \right .
  \end{equation*}
  with $(u^0,v^0)$ satisfying \eqref{cond_ini_propriedades} and $u^0_m \in C^{\infty}(\overline{\Omega})$ being mollifier regularizations of $u^0$ extended to $\mathbb{R}^N$ and having the following properties (see \cite{ViannaGuillen2023uniform} for more details):
  \begin{equation} \label{convergencia_dado_inicial_u_m}
    u^0_m\geq 0,\quad \int_\Omega u^0_m = \int_\Omega u^0,\quad 
    u^0_m \rightarrow u^0 \mbox{ strongly in } L^p(\Omega).
  \end{equation}
  
  
  \subsection{A \texorpdfstring{$\boldsymbol{L^{\infty}}$}{L\^{infty}} function bounding \texorpdfstring{$\boldsymbol{v_m}$}{v\_m} from above}
  
    In \cite{ViannaGuillen2023uniform}, where the uncontrolled model ($f \equiv 0$) is considered, a crucial step to prove the existence of a weak solution to \eqref{problema_P} as a limit of solutions of the truncated models \eqref{problema_P_m_controlado} is to get $m$-independent bounds for $\norma{v_m}{L^{\infty}(Q)}$. This fact also remains essential in our case, $f \not\equiv 0$. But, while in the case where $f \equiv 0$ this $m$-independent bound is obtained by straightforward calculations, it is not so obvious  now by considering a control $f$ with $f_+ \not\equiv 0$ in general.
    The next result will help us to build a function $w \in L^{\infty}$ bounding $v_m$ from above for all $m$. 
  \begin{lemma} \label{lema_problema_a_comparar}
    Let $\Omega$ be a bounded domain of $\mathbb{R}^3$ such that $\Gamma$ is of class $C^2$. Let $w^0 \in W^{2-2/q,q}(\Omega)$ and $\tilde{f} \in L^q(Q)$, for some $q > 5/2$. Then the problem
    \begin{equation} \label{problema_a_comparar}
      \left\{\begin{array}{l}
        \partial_t w - \Delta w = \tilde{f} \, w, \quad \mbox{ in } Q, \\
        \partial_\eta w |_{\Gamma}  = 0, \quad \mbox{on } (0,T) \times \Gamma, \\
        w(0,x)  = w^0 \  \mbox{in } \Omega,
      \end{array}\right.
    \end{equation}
    has a unique solution
    \begin{equation*}
      w \in C([0,T]; W^{2-2/q,q}(\Omega)) \cap L^q(0,T;W^{2,q}(\Omega)), \ \partial_t w \in L^q(Q),
    \end{equation*}
    In particular, there is a positive constant $C(\|\tilde{f}\|_{L^q(Q)}, \|w^0\|_{W^{2-2/q,q}(\Omega)})$ such that
    \begin{equation} \label{estimate_w_L_infty}
      \norma{w}{L^{\infty}(Q)} \leq C(\|\tilde{f}\|_{L^q(Q)}, \|w^0\|_{W^{2-2/q,q}(\Omega)}).
    \end{equation}
  \end{lemma}

\
  
  \noindent  In fact, inequality \eqref{estimate_w_L_infty} will provide an estimate for $\norma{v_m}{L^{\infty}(Q)}$ in terms of the control $f$. We remark that it is the main reason why we need to assume in this work that $f \in L^q(Q)$, for some $q > 5/2$.
   
  \begin{proof}[\bf Proof]
    The key idea here is the injection $W^{2-2/q,q}(\Omega) \subset L^{\infty}(\Omega)$, the reason why we suppose that $q > 5/2$. The proof is divided in two steps.
    
    \vspace{12pt}
    
    \noindent {\bf Step 1 (Existence and uniqueness of problem \eqref{problema_a_comparar}):}
    
    For any solution $w$ of \eqref{problema_a_comparar} such that
    \begin{equation} \label{strong_solution_problema_a_comparar}
      w \in L^{\infty}(0,T;H^1(\Omega)) \cap L^2(0,T;H^2(\Omega)), \ \partial_t w \in L^2(Q),
    \end{equation}
    we have
    \begin{equation} \label{L^p_estimates_w}
      \begin{array}{l}
      \displaystyle
        \norma{w(t)}{L^p(\Omega)}^p + \beta \int_0^t{\norma{\nabla [w(r)]^{p/2}}{L^2(\Omega)}^2 \ dr} \\[8pt]
        \displaystyle
        \leq C \norma{w^0}{L^p(\Omega)}^p \exp\left(C p^{5/2} \int_0^t \left(\|\tilde{f}(r)\|_{L^{5/2}(\Omega)}^{5/2} + 1\right) \ dr\right),
      \end{array}
    \end{equation}
    $a.e. \ t \in (0,T)$. In fact, we test the equation in \eqref{problema_a_comparar} by $p \, w^{p-1}$ and define $\tilde{w}: = w^{p/2}$. Using \eqref{interpolacao_norma_10/3}, we obtain
    \begin{align*}
      & \dfrac{d}{dt} \norma{\tilde{w}(t)}{L^2(\Omega)}^2 + \dfrac{4 p(p - 1)}{p^2} \int_0^t{\norma{\nabla \tilde{w}}{L^2(\Omega)}^2 \ dr} \leq p \int_{\Omega} \tilde{f} \ \tilde{w}^2 \ dx \\
      & \qquad \leq C\, p \,\|\tilde{f}\|_{L^{5/2}(\Omega)} \norma{\tilde{w}}{L^{10/3}(\Omega)}^2 \\
      & \qquad \leq C \,p\, \|\tilde{f}\|_{L^{5/2}(\Omega)} \norma{\tilde{w}}{L^2(\Omega)}^{4/5} \norma{\tilde{w}}{H^1(\Omega)}^{6/5} \\
      & \qquad \leq C(\delta) \, p^{5/2} \|\tilde{f}\|_{L^{5/2}(\Omega)}^{5/2} \norma{\tilde{w}}{L^2(\Omega)}^2 + \delta \norma{\tilde{w}}{L^2(\Omega)}^2 + \delta \norma{\nabla \tilde{w}}{L^2(\Omega)}^2.
    \end{align*}
    Hence, choosing $\delta > 0$ small enough to absorb the last term in the right hand side and going back to the notation $w$ we obtain
    \begin{equation*}
      \dfrac{d}{dt} \norma{w(t)}{L^p(\Omega)}^p + \beta \int_0^t{\norma{\nabla [w(r)]^{p/2}}{L^2(\Omega)}^2 \ dr} \leq C p^{5/2} ( \|\tilde{f}\|_{L^{5/2}(\Omega)}^{5/2} + 1) \norma{w(t)}{L^p(\Omega)}^p
    \end{equation*}
    and Gronwall's Lemma leads us to \eqref{L^p_estimates_w}. 
    
    For $\tilde{f}$ regular enough one can prove that \eqref{problema_a_comparar} has a unique solution satisfying \eqref{strong_solution_problema_a_comparar} by using Galerkin's method, for example. But accounting for the dependence of $w$ on $\|\tilde{f}\|_{L^{5/2}(\Omega)}^{5/2}$ given by \eqref{L^p_estimates_w}, we are actually able to prove that \eqref{problema_a_comparar} has a unique strong solution satisfying \eqref{strong_solution_problema_a_comparar} and \eqref{L^p_estimates_w} under a weaker assumption on the regularity of $\tilde{f}$. It is enough that $\tilde{f} \in L^{5/2}(Q)$, for instance. The uniqueness is proved by comparing two possibly distinct solutions of \eqref{problema_a_comparar} and concluding that they are in fact the same solution.
    
    \vspace{12pt}
    
    \noindent {\bf Step 2 (Proof of the $L^\infty$ estimate \eqref{estimate_w_L_infty}):}
    
    Since $\tilde{f} \in L^q(Q)$ and $q > 5/2$, \eqref{L^p_estimates_w} implies that there are $\tilde{q} \in (5/2,q)$ and a positive constant $\tilde{C}(\|\tilde{f}\|_{L^q(Q)}, \|w^0\|_{W^{2-2/q,q}(\Omega)})$ such that $f w \in L^{\tilde{q}}(Q)$ with
    \begin{equation} \label{norma_fw_dependencia_em_q}
      \|\tilde{f} \, w\|_{L^{\tilde{q}}(Q)} \leq \tilde{C}(\|\tilde{f}\|_{L^q(Q)}, \|w^0\|_{W^{2-2/q,q}(\Omega)}).
    \end{equation}
    From \eqref{norma_fw_dependencia_em_q} and \eqref{dependencia_continua_de_w_em_h} we can conclude,
    \begin{equation*}
      \norma{w}{C([0,T];W^{2-2/\tilde{q},\tilde{q}}(\Omega))} \leq C(\|\tilde{f}\|_{L^q(Q)}, \|w^0\|_{W^{2-2/q,q}(\Omega)}).
    \end{equation*}
    But since $\tilde{q} > 5/2$, we have $C([0,T];W^{2-2/\tilde{q},\tilde{q}}(\Omega)) \subset L^{\infty}(Q)$, hence
    \begin{equation*}
      \norma{w}{L^{\infty}(Q)} \leq C(\|\tilde{f}\|_{L^q(Q)}, \|w^0\|_{W^{2-2/q,q}(\Omega)}).
    \end{equation*}
    Finally, since $w \in L^{\infty}(Q)$, we have $\tilde{f} \ w \in L^q(Q)$ and therefore we can use Lemma \ref{lema_regularidade_eq_calor} to conclude that
    \begin{equation*}
      w \in C([0,T]; W^{2-2/q,q}(\Omega)) \cap L^q(0,T;W^{2,q}(\Omega)), \ \partial_t w \in L^q(Q).
    \end{equation*}
  \end{proof}

  
  \subsection{Existence for the controlled truncated problem and the first uniform estimates}
  
    \begin{theorem} \label{teo_existencia_P_m_controlado}
      Given $f_m$ and $(u^0_m,v^0)$ satisfying \eqref{properties_f_m}, \eqref{convergencia_dado_inicial_u_m} and $v^0 \in W^{2 - 2/q,q}(\Omega)$, respectively, there is a unique solution $(u_m,v_m)$ of \eqref{problema_P_m_controlado} with regularity
      \begin{equation*} 
        \begin{array}{c}
          u_m \in L^{\infty}(0,T;H^1(\Omega)) \cap L^2(0,T;H^2(\Omega)), \ \partial_t u_m \in L^2(0,T;L^2(\Omega)), \\
          v_m \in L^{\infty}(Q) \cap L^{\infty}(0,T;H^2(\Omega)),  \quad \Delta v_m \hbox{ and } \partial_t v_m \in L^2(0,T;H^1(\Omega)),
        \end{array}
      \end{equation*}
      and satisfying
      \begin{equation} \label{positivity_u_v}
        u_m(t,x), v_m(t,x) \geq 0, \ a.e. \ (t,x) \in Q,
      \end{equation}
      \begin{equation} \label{conservation_cells}
        \int_{\Omega}{u_m(t,x) \ dx} = \int_{\Omega}{u^0_m(x) \ dx}= \int_{\Omega}{u^0(x) \ dx}, \ a.e. \ t \in (0,T).
      \end{equation}
      Moreover, there is a positive, continuous and increasing function of $\norma{f}{L^q(Q)}$, $\mathcal{K}_1(\| f\|_{L^q(Q)}, \| v^0\|_{W^{2-2/q,q}(\Omega)})$, also independent of $m$, such that 
      \begin{equation} \label{first_estimates_v_m}
        \norma{v_m}{L^{\infty}(Q)}, \norma{v_m}{L^2(0,T; H^1(\Omega))} \leq \mathcal{K}_1(\norma{f}{L^q(Q)}, \norma{v^0}{W^{2-2/q,q}(\Omega)}).
      \end{equation}
    \end{theorem}
    \begin{proof}[\bf Proof]
      Concerning the proof of existence and uniqueness of solution to \eqref{problema_P_m_controlado}, the truncation $T^m(\cdot)$ simplifies the treatment of the chemotaxis and consumption terms, $- \nabla \cdot (T^m(u_m) \nabla v_m)$ and $- T^m(u_m)^s v_m$ respectively, and one can deal with the control term $f_m v_m$, likewise in the proof of existence of problem  \eqref{problema_a_comparar}, in Lemma \ref{lema_problema_a_comparar}. The uniqueness is proved by comparing two possibly different solutions. The regularity of the solution pair $(u_m,v_m)$ and properties \eqref{positivity_u_v} and \eqref{conservation_cells} can be proved following the ideas in \cite{ViannaGuillen2023uniform}. Finally, we prove \eqref{first_estimates_v_m}, beginning by the estimate in the $L^{\infty}$-norm. Using the already proved property \eqref{positivity_u_v} of $v_m$ in the $v_m$-equation of \eqref{problema_P_m_controlado}, we obtain
      \begin{equation} \label{desigualdade_v_m_f_+}
        \partial_t v_m - \Delta v_m \leq (f_m)_+ \, v_m \quad a.e. \ (t,x) \in Q.
      \end{equation}
      On the other hand, accounting for Lemma \ref{lema_problema_a_comparar}, as well as \eqref{problema_a_comparar}, with $\tilde{f} = f_m)_+$ and $w^0 = v^0$, we consider $w$ satisfying
      \begin{equation} \label{problema_a_comparar_com_tilde_f_igual_a_f_+}
        \partial_t w - \Delta w = (f_m)_+ \, w \quad a.e. \ (t,x) \in Q,
      \end{equation}
      with $\partial_\eta w  |_{\Gamma} = 0$ and $w(0,x) = v^0(x)$. Subtracting \eqref{problema_a_comparar_com_tilde_f_igual_a_f_+} from \eqref{desigualdade_v_m_f_+} we conclude that $(v_m - w)$ satisfies
      \begin{equation*}
        \left \{ \begin{array}{l}
          \partial_t (v_m - w) - \Delta (v_m - w) \leq 
          (f_m)_+ \, (v_m - w) \quad a.e. \ (t,x) \in Q, \\[6pt]
          \partial_\eta (v_m - w)  |_{\Gamma} = 0, \quad 
          (v_m - w)(0,x) = 0.
        \end{array} \right.
      \end{equation*}
      Multiplying the above inequality by $(v_m - w)_+$ and using \eqref{interpolacao_norma_10/3} to estimate the right hand side term leads us to $(v_m - w)_+(t,x) = 0$ $a.e. \ (t,x) \in Q$, that is, $v_m(t,x) \leq w(t,x)$ $a.e. \ (t,x) \in Q$. Then, taking the $L^q$-norm inequality of \eqref{properties_f_m} into account, the bound in the $L^{\infty}$-norm for $v_m$ is a consequence of the estimate for $w$ given in \eqref{estimate_w_L_infty}, with $\tilde{f} = (f_m)_+$ and $w^0 = v^0$, from Lemma \ref{lema_problema_a_comparar}. The bound of $v_m$ in the norm of $L^2(0,T;H^1(\Omega))$ is obtained by testing the $v_m$-equation of \eqref{problema_P_m_controlado} by $v_m$ and conveniently estimating the term on the right hand side using Holder's inequality, the interpolation inequality \eqref{interpolacao_norma_10/3} and Young's inequality.
    \end{proof}


  \subsection{Energy inequality}
    
    Analogously to \cite{ViannaGuillen2023uniform}, we consider the variable $z_m(t,x) = \sqrt{v_m(t,x) + \alpha^2}$ and the rewritten problem 
    {\small \begin{equation} \label{problema_P_u_m_z_m_controlado}
      \begin{array}{rl}
        \partial_t u_m - \Delta u_m & = - \nabla \cdot (T^{m}(u_m) \nabla (z_m)^2) \\[6pt]
        \partial_t z_m - \Delta z_m - \dfrac{\norm{\nabla z_m}{}^2}{z_m} &
         = - \dfrac{1}{2} T^{m}(u_m)^s \left(z_m - \dfrac{\alpha^2}{z_m} \right) + \dfrac{1}{2} f_m \left(z_m - \dfrac{\alpha^2}{z_m} \right) 1_{\Omega_c} \\[6pt]
        \partial_\eta u_m  |_{\Gamma} & =  \partial_\eta z_m  |_{\Gamma} = 0 \\[6pt]
        u_m(0) & = u^0_m, \quad z_m(0) = \sqrt{v^0 + \alpha^2},
      \end{array}
    \end{equation}}
    which is equivalent to the controlled truncated problem \eqref{problema_P_m_controlado}. From the equivalence of \eqref{problema_P_m_controlado} and \eqref{problema_P_u_m_z_m_controlado} and from the results given in Lemma \ref{teo_existencia_P_m_controlado}, we have the following. 
    \begin{corollary} \label{coro_existencia_P_u_m_z_m_controlado}
      Given $f_m$ and $(u^0_m,v^0)$ satisfying \eqref{properties_f_m}, \eqref{convergencia_dado_inicial_u_m} and $v^0 \in W^{2 - 2/q,q}(\Omega)$, respectively, there is a unique solution $(u_m,z_m)$ of \eqref{problema_P_u_m_z_m_controlado} with regularity
      \begin{equation*} 
        \begin{array}{c}
          u_m \in L^{\infty}(0,T;H^1(\Omega)) \cap L^2(0,T;H^2(\Omega)), \ \partial_t u_m \in L^2(0,T;L^2(\Omega))), \\
          z_m \in L^{\infty}(Q) \cap L^{\infty}(0,T;H^2(\Omega)), \quad \Delta z_m \hbox{ and } \partial_t z_m \in L^2(0,T;H^1(\Omega)),
        \end{array}
      \end{equation*}
      and satisfying the $m$-uniform estimates
      \begin{equation*}
        u_m(t,x) \geq 0 \ \mbox{ and } \ z_m(t,x) \geq \alpha, \ a.e. \ (t,x) \in Q,
      \end{equation*}
      \begin{equation*}
        \int_{\Omega}{u_m(t,x) \ dx} = \int_{\Omega}{u^0_m(x) \ dx} = \int_{\Omega}{u^0(x) \ dx}, \ a.e. \ t \in (0,T),
      \end{equation*}
      \begin{equation} \label{first_estimates_z_m}
        \norma{z_m}{L^{\infty}(Q)}, \norma{z_m}{L^2(0,T;H^1(\Omega))} \leq \mathcal{K}_1(\norma{f}{L^q(Q)}, \norma{v^0}{W^{2-2/q,q}(\Omega)}).
      \end{equation}
    \end{corollary}

\

    Using this change of variables, we obtain an energy inequality involving the control $f$. In this subsection, in order to simplify the notation, we drop the $m$ subscript and denote the solution $(u_m,z_m)$ of \eqref{problema_P_u_m_z_m_controlado} by $(u,z)$. We begin with the following lemma.
    \begin{lemma} \label{lema_tratamento_equacao_z}
      The solution $(u,z)$ of \eqref{problema_P_u_m_z_m_controlado}, satisfies the inequality
      \begin{equation*}
        \begin{array}{l}
          \dfrac{1}{2} \dfrac{d}{dt} \norma{\nabla z}{L^2(\Omega)}^2 + C_1 \Big ( \D{\int_{\Omega}}{\norm{D^2 z}{}^2 \ dx} + \D{\int_{\Omega}}{\frac{\norm{\nabla z}{}^4}{z^2} \ dx} \Big ) + \dfrac{1}{2} \D{\int_{\Omega}{T^{m}(u)^s \norm{\nabla z}{}^2 \ dx}} \\
          \qquad \leq \dfrac{s}{2} \alpha \D{\int_{\Omega}}{T^{m}(u)^{s-1} \norm{\nabla z}{} \norm{\nabla T^{m}(u)}{} \ dx} + \mathcal{K}_1^2 \norma{f}{L^2(\Omega)}^2 \\
          \qquad \qquad + C_2 \norma{\nabla z}{L^2(\Omega)}^2 + \dfrac{s}{4} \D{\int_{\Omega}}{T^{m}(u)^{s-1} \nabla (z^2) \cdot \nabla T^{m}(u) \ dx}.
        \end{array}
      \end{equation*}
    \end{lemma}
    \begin{proof}[\bf Proof]
      The proof of this lemma is analogous to the proof of Lemma 23 of \cite{ViannaGuillen2023uniform}. We also use property \eqref{properties_f_m} which says that $\norma{f_m}{L^2(Q)} \leq \norma{f}{L^2(Q)}$.
    \end{proof}

    Now we need to prove $m$-independent estimates for $(u_m,z_m)$. These estimates will be obtained from an energy inequality as in \cite[Subsection 5.1 and 5.2]{ViannaGuillen2023uniform}, with the modifications due to the treatment of the term related to the control $f$. We remark that we will use property \eqref{properties_f_m} to keep the dependence on the norms of the control $f$ rather than in terms of the norms of the mollified functions $f_m$. As in \cite{ViannaGuillen2023uniform}, we have to separate the proof of the energy inequality for the cases $s = 1$, $s \in (1,2)$ and $s \geq 2$.  In fact, the test function used for the $u$-equation is $g'(u)$, where $g''(u) = u^{s-2}$. Note that, if $s = 1$, both $g'(u)$ and $g''(u)$ are singular at $u = 0$; if $s \in (1,2)$, then only $g''(u)$ is singular at $u = 0$; and if $s \geq 2$ then neither $g'(u)$ nor $g''(u)$ have any singularity.
    
    Next we consider the function $g_m$ defined by $g_m(r) = \D{\int_0^r}{g'_m(\theta) \ d\theta}$, where $g_m'(\theta)$ is defined for $\theta \geq 0$ by
    \begin{equation*}
      g_m'(\theta) = \left \{ 
      \begin{array}{rl}
        ln(T^m(\theta) + 1), & \mbox{if } s = 1, \\
        \dfrac{T^m(\theta)^{s-1}}{(s-1)}, & \mbox{if } s > 1,
      \end{array}
      \right.
    \end{equation*}
    and the energy
    \begin{equation} \label{energia_E_m}
      E_m(u,z)(t) = \frac{s}{4} \int_{\Omega}{g_m(u(t,x)) \ dx} + \frac{1}{2} \int_{\Omega}{\norm{\nabla z(t,x)}{}^2} \ dx.
    \end{equation}
    
    \begin{lemma}[\bf Energy inequality for $\boldsymbol{s = 1}$] \label{lemma_estimativa_u_v_m_1_s=1}
        The solution $(u, z)$ of the problem \eqref{problema_P_u_m_z_m_controlado} satisfies, for sufficiently small $\alpha > 0$,
        \begin{equation} \label{estimativa_u_v_m_1_s=1}
          \begin{array}{l}
            \dfrac{d}{dt} E_m(u,z)(t) + \beta \D{\int_{\Omega}}{\norm{\nabla [T^{m}(u) + 1]^{1/2}}{}^2 \ dx} + \dfrac{1}{4} \D{\int_{\Omega}}{T^{m}(u) \norm{\nabla z}{}^2 \ dx} \\
            + \beta \Big ( \D{\int_{\Omega}}{\norm{D^2 z}{}^2 \ dx} + \D{\int_{\Omega}}{\frac{\norm{\nabla z}{}^4}{z^2} \ dx} \Big ) \leq C (\mathcal{K}^2) \norma{\nabla z}{L^2(\Omega)}^2 + \mathcal{K}_1^2 \norma{f}{L^2(\Omega)}^2.
          \end{array}
        \end{equation}
      \end{lemma}
      \begin{proof}[\bf Proof]
        We follow \cite{ViannaGuillen2023uniform}, pointing out the most relevant steps to deal with the control term and make explicit the dependence on  $\mathcal{K}_1 = \mathcal{K}_1(\| f\|_{L^q(Q)}, \|v^0\|_{W^{2-2/q,q}(\Omega)})$, from Corollary \ref{coro_existencia_P_u_m_z_m_controlado}. By testing the $u$-equation of problem \eqref{problema_P_u_m_z_m_controlado} by $ ln(T^{m}(u) + 1)$, we obtain
        \begin{equation*}
          \frac{d}{dt} \int_{\Omega}{g_m(u) \ dx} + \int_{\Omega}{\frac{(T^{m})'(u)}{T^{m}(u) + 1} \norm{\nabla u}{}^2 \ dx} = \prodl{\frac{T^{m}(u)}{T^{m}(u) + 1} \nabla (z^2)}{\nabla T^{m}(u)}.
        \end{equation*}
        Since $0 \leq (T^{m})'(u) \leq C$, we have $((T^{m})'(u))^2 \leq C (T^{m})'(u)$, and we can write
        \begin{align*}
          \int_{\Omega}{\frac{(T^{m})'(u)}{T^{m}(u) + 1} \norm{\nabla u}{}^2 \ dx} \geq C \int_{\Omega}{\frac{((T^{m})'(u))^2}{T^{m}(u) + 1} \norm{\nabla u}{}^2 \ dx} \geq C \int_{\Omega}{\norm{\nabla [T^{m}(u) + 1]^{1/2}}{}^2 \ dx}.
        \end{align*}
        
        Hence, using that $\dfrac{1}{(T^{m}(u) + 1)} \leq \dfrac{1}{\sqrt{T^{m}(u) + 1}}$, we have
        \begin{align*}
          & \frac{d}{dt} \int_{\Omega}{g_m(u) \ dx} + C \int_{\Omega}{\norm{\nabla [T^{m}(u) + 1]^{1/2}}{}^2 \ dx} = 2 \prodl{\frac{T^{m}(u) + 1 - 1}{T^{m}(u) + 1} z \nabla z}{\nabla T^{m}(u)} \\
          & = \prodl{\nabla (z^2)}{\nabla T^{m}(u)} - 2 \prodl{z \nabla z}{\frac{\nabla T^{m}(u)}{T^{m}(u) + 1}} \\
          & \leq \prodl{\nabla (z^2)}{\nabla T^{m}(u)} + 2 \norma{z}{L^{\infty}(\Omega)} \norma{\nabla z}{L^2(\Omega)} \norma{\nabla [T^{m}(u) + 1]^{1/2}}{L^2(\Omega)}.
        \end{align*}
        Using Young's inequality and \eqref{first_estimates_z_m}, we arrive at
        \begin{equation} \label{estimativa_u_m_1_s=1}
          \begin{array}{rl}
            & \dfrac{d}{dt} \D{\int_{\Omega}}{g_m(u) \ dx} + C \D{\int_{\Omega}}{\norm{\nabla [T^{m}(u) + 1]^{1/2}}{}^2 \ dx} \\[6pt]
            &\qquad  \leq \prodl{\nabla (z^2)}{\nabla T^{m}(u)} + \mathcal{K}_1^2 \norma{\nabla z}{L^2(\Omega)}^2.
          \end{array}
        \end{equation}
        
        If we add the inequality of Lemma \ref{lema_tratamento_equacao_z}, for $s = 1$, to $1/4$ times \eqref{estimativa_u_m_1_s=1}, then the terms $\D{\int_{\Omega}{\nabla T^{m}(u) \cdot \nabla (z^2) \ dx}}$ cancel and we obtain
        \begin{align}
          & \nonumber \frac{d}{dt} \Big [ \frac{1}{4} \int_{\Omega}{g_m(u) \ dx} + \frac{1}{2} \norma{\nabla z}{L^2(\Omega)}^2 \Big ] + C \int_{\Omega}{\norm{\nabla [T^{m}(u) + 1]^{1/2}}{}^2 \ dx} \\
          & \nonumber \qquad \qquad + \frac{1}{2} \int_{\Omega}{T^{m}(u) \norm{\nabla z}{}^2 \ dx} + C_1 \Big ( \int_{\Omega}{\norm{D^2 z}{}^2 \ dx} + \int_{\Omega}{\frac{\norm{\nabla z}{}^4}{z^2} \ dx} \Big ) \\
          & \qquad \leq \frac{\sqrt{\alpha}}{2} \int_{\Omega}{\norm{\nabla z}{} \norm{\nabla T^{m}(u)}{} \ dx} + \mathcal{K}_1^2 \norma{f}{L^2(\Omega)}^2 + (C_2 + \mathcal{K}_1^2) \norma{\nabla z}{L^2(\Omega)}^2 \label{aux_inequality_u_z_s=1} \\
          & \nonumber \qquad \leq \int_{\Omega}{\alpha \norm{\nabla [T^{m}(u) + 1]^{1/2}}{} \norm{\sqrt{T^{m}(u) + 1}}{} \norm{\nabla z}{} \ dx} \\
          & \nonumber \qquad \qquad + \mathcal{K}_1^2 \norma{f}{L^2(\Omega)}^2 + (C_2 + \mathcal{K}_1^2) \norma{\nabla z}{L^2(\Omega)}^2.
        \end{align}
        We can deal with the first term in the right hand side of the inequality using Holder's and Young's inequality,
        \begin{align*}
          & \int_{\Omega}{\alpha \norm{\nabla [T^{m}(u) + 1]^{1/2}}{} \norm{\sqrt{T^{m}(u) + 1}}{} \norm{\nabla z}{} \ dx} \leq 
          \alpha^2 \ C(\delta) \int_{\Omega}{T^{m}(u) \norm{\nabla z}{}^2 \ dx} \\
          & \qquad + \delta \norma{\nabla [T^{m}(u) + 1]^{1/2}}{L^2(\Omega)}^2 + \alpha^2 \ C(\delta) \int_{\Omega}{\norm{\nabla z}{}^2 \ dx}.
        \end{align*}
        
        Therefore, we can first choose $\delta > 0$ and then $\alpha > 0$ sufficiently small in order to use the terms on the left hand side of inequality \eqref{aux_inequality_u_z_s=1} to absorb  the first two terms on the right hand side of the above inequality and finally obtain the desired inequality \eqref{estimativa_u_v_m_1_s=1}.
      \end{proof}
      
      Now we  obtain the energy inequalities for $s \in (1,2)$ and for $s \geq 2$, respectively. Analogously to Lemma \ref{lemma_estimativa_u_v_m_1_s=1}, we follow the ideas of \cite{ViannaGuillen2023uniform}, making the necessary changes in order to deal with the control term and to make explicit the dependence on the positive constant $\mathcal{K}_1 = \mathcal{K}_1(\| f\|_{L^q(Q)}, \| v^0\|_{W^{2-2/q,q}(\Omega)})$, from Corollary \ref{coro_existencia_P_u_m_z_m_controlado}. Since these changes were covered in Lemma \ref{lemma_estimativa_u_v_m_1_s=1}, next we will state the results, skipping their proofs. 
      
      \begin{lemma}[\bf Energy inequality for $\boldsymbol{s \in (1,2)}$]
        The solution $(u, z)$ of the problem \eqref{problema_P_u_m_z_m_controlado} satisfies, for sufficiently small $\alpha > 0$,
        \begin{equation}
          \begin{array}{l}
            \dfrac{d}{dt} E_m(u,z)(t) + \beta \D{\int_{\Omega}}{\norm{\nabla [T^{m}(u) + 1]^{s/2}}{}^2 \ dx} + \dfrac{1}{4} \D{\int_{\Omega}}{T^m(u)^s \norm{\nabla z}{}^2 \ dx} \\[6pt]
            + \beta \Big ( \D{\int_{\Omega}}{\norm{D^2 z}{}^2 \ dx} + \D{\int_{\Omega}}{\frac{\norm{\nabla z}{}^4}{z^2} \ dx} \Big ) \leq C (\mathcal{K}_1^2) \norma{\nabla z}{L^2(\Omega)}^2 + \mathcal{K}_1^2 \norma{f}{L^2(\Omega)}^2.
          \end{array}
          \label{estimativa_u_v_m_1_s_intermediario}
        \end{equation}
        \label{lemma_estimativa_u_v_m_1_s_intermediario}
      \end{lemma}
      \begin{remark}
        In \cite{ViannaGuillen2023uniform}, in the lemma where the authors prove the energy inequality for $s \in (1,2)$, the term $\D{\int_{\Omega}}{\norm{\nabla [T^{m}(u) + 1/j]^{s/2}}{}^2 dx}$ is estimated by 
        \begin{equation*}
          \D{\int_{\Omega}}{\norm{\nabla [T^{m}(u) + 1/j]^{s/2}}{}^2 dx} \geq 0, \mbox{ for all } j \in \mathbb{N},
        \end{equation*}
        but it can be estimated by
        \begin{equation*}
          \D{\int_{\Omega}}{\norm{\nabla [T^{m}(u) + 1/j]^{s/2}}{}^2 dx} \geq \D{\int_{\Omega}}{\norm{\nabla [T^{m}(u) + 1]^{s/2}}{}^2 dx}, \mbox{ for all } j \in \mathbb{N},
        \end{equation*}
        instead, yielding \eqref{estimativa_u_v_m_1_s_intermediario}.
        \hfill $\square$
      \end{remark}
      
      \begin{lemma}[\bf Energy inequality for $\boldsymbol{s \geq 2}$]
        The solution $(u, z)$ of the problem \eqref{problema_P_u_m_z_m_controlado} satisfies, for sufficiently small $\alpha > 0$,
        \begin{equation}
          \begin{array}{l}
            \dfrac{d}{dt} E_m(u,z)(t) + \D{\int_{\Omega}}{\norm{\nabla [T^m(u)]^{s/2}}{}^2  dx} + \dfrac{1}{4} \D{\int_{\Omega}}{T^m(u)^s \norm{\nabla z}{}^2 dx} \\
            + \beta \Big ( \D{\int_{\Omega}}{\norm{D^2 z}{}^2  dx} + \D{\int_{\Omega}}{\frac{\norm{\nabla z}{}^4}{z^2}  dx} \Big ) \leq C (\mathcal{K}_1^2) \norma{\nabla z}{L^2(\Omega)}^2 + \mathcal{K}_1^2 \norma{f}{L^2(\Omega)}^2.
          \end{array}
          \label{estimativa_u_v_m_1_s_geq_2}
        \end{equation}
        \label{lemma_estimativa_u_v_m_1_s_geq_2}
      \end{lemma}


  \subsection{\texorpdfstring{$\boldsymbol{m}$}{m}-independent estimates and passage to the limit as \texorpdfstring{$\boldsymbol{m \to \infty}$}{m goes to infinity}}
  \label{subsec:estimativas_u_m_v_m_e_passagem_ao_limite}
    
    In the present subsection we go back to the notation $(u_m,z_m)$ and $(u_m,v_m)$ to the solution of \eqref{problema_P_u_m_z_m_controlado} and \eqref{problema_P_m_controlado}, respectively.


  \subsubsection{\bf \texorpdfstring{$\boldsymbol{m}$}{m}-independent estimates for \texorpdfstring{$\boldsymbol{\nabla v_m}$}{grad(v\_m)}}
  \label{subsec:estimativas_v_m}
    
    We will integrate the energy inequalities \eqref{estimativa_u_v_m_1_s=1}, \eqref{estimativa_u_v_m_1_s_intermediario} and \eqref{estimativa_u_v_m_1_s_geq_2} with respect to $t$, from $0$ to $T > 0$. We take into account that, because of Corollary \ref{coro_existencia_P_u_m_z_m_controlado}, we have the following bounds independently of  $m$:
    \begin{equation*}
      \nabla z_m \mbox{ is bounded in } L^{2}(Q)
    \end{equation*}
    and
    \begin{equation} \label{limitacao_por_baixo_e_por_cima_z_m}
      0 < \alpha \leq z_m(t,x) \leq \mathcal{K}_1, \ a.e. \ (t,x) \in Q.
    \end{equation}
    We also use the hypothesis \eqref{cond_ini_propriedades} on the initial data $u^0,v^0$ to prove that the energy given in \eqref{energia_E_m} at time $t = 0$, $E_m(u_m,z_m)(0)$, is also bounded, independently of $m$. Thus we conclude that
      \begin{equation*}
        \nabla z_m \mbox{ is bounded in } L^{\infty}(0,T;L^2(\Omega)) \cap L^4(Q),
      \end{equation*}
      \begin{equation*}
        T^m(u_m)^{s/2} \nabla z_m \mbox{ and } \Delta z_m \mbox{ are bounded in } L^2(Q).
      \end{equation*}
      But using the fact that $z_m = \sqrt{v_m + \alpha^2}$ and \eqref{limitacao_por_baixo_e_por_cima_z_m} we can conclude that
      \begin{equation}
        \nabla v_m \mbox{ is bounded in } L^{\infty}(0,T;L^2(\Omega)) \cap L^4(Q),
      \label{limitacao_aux_Dv_m_s_intermediario}
      \end{equation}
      \begin{equation}
        T^m(u_m)^{s/2} \nabla v_m \mbox{ and } \Delta v_m \mbox{ are bounded in } L^2(Q).
        \label{limitacao_aux_Delta_v_m_s_intermediario}
      \end{equation}


  \subsubsection{\bf Case \texorpdfstring{$\boldsymbol{s \in [1,2)}$}{s in [1,2)}}
    \label{subsec:estimativas_u_m_s<2}
    
    First, following \cite{ViannaGuillen2023uniform}, to which we refer the reader that might be interested in more details, we prove the existence of weak solution $(u,v)$ to \eqref{problema_P_controlado}. Next, to conclude the proof of Theorem \ref{teo_existencia_solucao_fraca}, letting $z = \sqrt{v + \alpha^2}$, we prove the energy inequality \eqref{desigualdade_integral_limite}.
  

\

\noindent {\bf Existence of weak solution to \eqref{problema_P_controlado}:}
  
    In order to prove the existence of a weak solution $(u,v)$ to \eqref{problema_P_controlado}, first we obtain $m$-independent estimates for the solutions $(u_m,v_m)$ of \eqref{problema_P_m_controlado} and then we use compactness results in weak*, weak and strong topologies to pass to the limit as $m \to \infty$.
    
    Let
      \begin{equation*}
        \qquad g'(r) =  \left \{
        \begin{array}{cc}
          ln(r) & \mbox{ if } s = 1, \\
          r^{s-1}/(s-1) & \mbox{ if } s \in (1,2),
        \end{array}
        \right .
         \forall r > 0.
      \end{equation*}
      and let
      \begin{equation*}
        g(r) =  \int_0^r{g'(\theta) \ d\theta} = \left \{
        \begin{array}{cc}
          r ln(r) - (r - 1) & \mbox{ if } s = 1, \\
          r^s/s(s-1) & \mbox{ if } s \in (1,2).
        \end{array}
        \right .
      \end{equation*}
      Notice that $g''(r) = r^{s-2}, \ \forall r > 0$, in all cases. 
      We test the $u_m$-equation of \eqref{problema_P_m_controlado} by $g'(u_m + 1)$ and obtain
      \begin{align*}
        \dfrac{d}{dt} \int_{\Omega} g(u_m + 1) \ dx & + \dfrac{4}{s^2} \int_{\Omega}{\norm{\nabla [u_m + 1]^{s/2}}{}^2 \ dx} \\
        & = \int_{\Omega}{T^m(u_m) (u_m + 1)^{s/2-1} \nabla v_m  \cdot \nabla u_m\,  (u_m + 1)^{s/2-1} \ dx} \\
        & = \dfrac{2}{s} \int_{\Omega}{\dfrac{T^m(u_m)^{1-s/2}}{(u_m + 1)^{1-s/2}} T^m(u_m)^{s/2} \nabla v_m \cdot \nabla [u_m + 1]^{s/2} \ dx} \\
        & \leq \dfrac{2}{s} \Big ( \int_{\Omega}{T^m(u_m)^s \norm{\nabla v_m}{}^2 \ dx} \Big )^{1/2} \Big ( \int_{\Omega}{\norm{\nabla [u_m + 1]^{s/2}}{}^2 \ dx} \Big )^{1/2}
      \end{align*}
      and thus we have
 $$
        \dfrac{d}{dt} \int_{\Omega}{g(u_m + 1) dx}  + \dfrac{2}{s^2} \int_{\Omega}{\norm{\nabla [u_m + 1]^{s/2}}{}^2 dx}
        \leq \dfrac{1}{4} \int_{\Omega}{T^m(u_m)^s \norm{\nabla v_m}{}^2 dx}.
 $$
      Integrating with respect to $t$ from $0$ to $T$ we obtain
      \begin{align*}
        \int_{\Omega}{g(u_m(T) + 1) dx} + \dfrac{2}{s^2} \int_{0}^T{\int_{\Omega}{\norm{\nabla [u_m + 1]^{s/2}}{}^2 dx} dt} \\
        \leq \dfrac{1}{4} \int_0^T{\int_{\Omega}{T^m(u_m)^s \norm{\nabla v_m}{}^2 dx}  dt} + \int_{\Omega}{g(u^0 + 1)  dx}.
      \end{align*}
      
      Then, because of \eqref{first_estimates_v_m}, \eqref{cond_ini_propriedades} and the definition of $g$ and \eqref{limitacao_aux_Delta_v_m_s_intermediario} we conclude that
      \begin{equation}
        (u_m + 1)^{s/2} \mbox{ is bounded in } L^{\infty}(0,T;L^2(\Omega)) \cap L^2(0,T;H^1(\Omega)).
        \label{limitacao_u_m_L^s}
      \end{equation}
      Using the Sobolev inequality $H^1(\Omega) \subset L^6(\Omega)$ and interpolation inequalities we obtain
      \begin{equation*}
        (u_m)^{s/2} \mbox{ is bounded in } L^{10/3}(Q).
      \end{equation*}
      The latter and \eqref{limitacao_u_m_L^s} imply that
      \begin{equation}
        u_m \mbox{ is bounded in } L^{\infty}(0,T;L^s(\Omega)) \cap L^{5s/3}(Q).
        \label{limitacao_u_m_L_5s/3}
      \end{equation}
      
      From \eqref{limitacao_u_m_L_5s/3} we can conclude, using the $v_m$-equation of \eqref{problema_P_m_controlado} that
      \begin{equation*} 
        \partial_t v_m \mbox{ is bounded in } L^{5/3}(Q).
      \end{equation*}
      
      Reminding that $s \in [1,2)$, if we use \eqref{limitacao_u_m_L^s} and \eqref{limitacao_u_m_L_5s/3} in the relation
      \begin{equation*} 
        \nabla u_m = \nabla (u_m + 1) = \nabla \big ( (u_m + 1)^{s/2} \big )^{2/s} = \dfrac{2}{s} (u_m + 1)^{1 - s/2} \ \nabla (u_m + 1)^{s/2}.
      \end{equation*}
      then we also have
      \begin{equation}
        u_m \mbox{ is bounded in } L^{5s/(3+s)}(0,T;W^{1,5s/(3+s)}(\Omega)).
        \label{limitacao_W_5s_3+s}
      \end{equation}
      
      Considering the chemotaxis term of the $u_m$-equation of \eqref{problema_P_m_controlado}, we can write $T^m(u_m) \nabla v_m$ as
  \begin{equation*}
    T^m(u_m) \nabla v_m = T^m(u_m)^{1-s/2} T^m(u_m)^{s/2} \nabla v_m.
  \end{equation*}
  Then, we have $T^m(u_m)^{1-s/2}$ bounded in $L^{10s/(6-3s)}(Q)$, because of \eqref{limitacao_u_m_L^s}, and $T^m(u_m)^{s/2} \nabla v_m$ bounded in $L^2(Q)$, because of \eqref{limitacao_aux_Delta_v_m_s_intermediario}, and hence we can conclude that
  \begin{equation}
    T^m(u_m) \nabla v_m \mbox{ is bounded in } L^{5s/(3+s)}(Q).
    \label{limitacao_termo_chemotaxis_s}
  \end{equation}
  Then, if we consider the $u_m$-equation of \eqref{problema_P_m_controlado}, from \eqref{limitacao_W_5s_3+s} and \eqref{limitacao_termo_chemotaxis_s} we conclude that 
  \begin{equation*} 
    \partial_t u_m \mbox{ is bounded in } L^{5s/(3+s)}(0,T;(W^{1,5s/(4s-3)}(\Omega))').
  \end{equation*}
  
      Now we are going to obtain compactness for $\{ u_m \}$ which is necessary in order to pass to the limit as $m \to \infty$ in the nonlinear terms of the equations of \eqref{problema_P_m_controlado}. 
  
  We observe that $W^{1,5s/(3+s)}(\Omega) \subset L^q(\Omega)$, with continuous embedding for $q = 15s/(9-2s)$ and compact embedding for $q \in [1,15s/(9-2s))$. Then, since $s \in [1,2)$, we have $5s/3 < 15s/(9-2s)$ and therefore the embedding $W^{1,5s/(3+s)}(\Omega) \subset L^{5s/3}(\Omega)$ is compact. Note also that $q = 5s/3 \geq 5/3 > 1$.
  
  Now we can use Lemma \ref{lema_Simon} with 
 \[ X = W^{1,5s/(3+s)}(\Omega),\quad
      B = L^{5s/3}(\Omega),  \quad Y = \big ( W^{1,5s/(4s-3)}(\Omega) \big )'\]
  and $q = 5s/3$, to conclude that there is a subsequence of $\{ u_m \}$ (still denoted by $\{ u_m \}$) and a limit function $u$ such that
  \begin{equation*}
    u_m \longrightarrow u \mbox{ weakly in } L^{5s/(3+s)}(0,T;W^{1,5s/(3+s)}(\Omega)),
  \end{equation*}
  and
  \begin{equation} \label{convergencia_u_m_s_intermediario}
    u_m \longrightarrow u \mbox{ strongly in } L^p(Q), \ \forall p \in [1,5s/3).
  \end{equation}
  
  Using the Dominated Convergence Theorem we can conclude from \eqref{convergencia_u_m_s_intermediario} that
  \begin{equation} \label{convergencia_T^m_u_s_intermediario}
    T^m(u_m) \rightarrow u \mbox{ strongly in } L^p(Q), \ \forall p \in [1,5s/3).
  \end{equation}
  It stems from the convergence \eqref{convergencia_T^m_u_s_intermediario} and Lemma \ref{lema_convergencia_w_elevado_a_s} that
  \begin{equation} \label{convergencia_T^m_u_elevado_a_s_s_intermediario}
    (T^m(u_m))^s \rightarrow u^s \mbox{ strongly in } L^q(Q), \ \forall q \in [1,5/3).
  \end{equation}
  
  The convergence of $v_m$ is better. There is a subsequence of $\{ v_m \}$ (still denoted by $\{ v_m \}$) and a limit function $v$ such that
  \begin{equation}
    \begin{array}{c}
      v_m \rightarrow v \mbox{ weakly* in } L^{\infty}(Q) \cap L^{\infty}(0,T;H^1(\Omega)), \\
      v_m \rightarrow v \mbox{ weakly in } L^2(0,T;H^2(\Omega)), \\
      \nabla v_m \rightarrow \nabla v \mbox{ weakly in } L^4(Q), \\
      \mbox{and } \partial_t v_m \rightarrow \partial_t v \mbox{ weakly in } L^{5/3}(Q).
    \end{array}
    \label{convergencia_v_m_s_intermediario}
  \end{equation}
  
  Now we are going to use the weak and strong convergences obtained so far to pass to the limit as $m \to \infty$ in the equations of problem \eqref{problema_P_m_controlado}. Since passing to the limit in the linear terms is simpler, as well as in the term $f_m v_m$, because of the strong convergence of $f_m$ given in \eqref{properties_f_m}, we focus on the nonlinear terms of the equations. We are going to identify the limits of the nonlinear terms related to chemotaxis and consumption,
  \begin{equation*}
    T^m(u_m) \nabla v_m \mbox{ and } T^m(u_m)^s v_m,
  \end{equation*}
  respectively, with 
  \begin{equation*}
    u \nabla v \mbox{ and } u^s v.
  \end{equation*}
  In fact, considering the chemotaxis term, because of \eqref{convergencia_T^m_u_s_intermediario}, \eqref{limitacao_aux_Dv_m_s_intermediario} and \eqref{convergencia_v_m_s_intermediario}, we can conclude that
  \begin{equation*}
    T^m(u_m) \nabla v_m \longrightarrow u \nabla v \mbox{ weakly in } L^{20s/(5s+12)}(Q).
  \end{equation*}
  
  Considering now the consumption term, considering \eqref{convergencia_T^m_u_elevado_a_s_s_intermediario} and \eqref{convergencia_v_m_s_intermediario} we conclude that
  \begin{equation*}
    T^m(u_m)^s v_m \longrightarrow u^s v \mbox{ weakly in } L^{5/3}(Q).
  \end{equation*}
  
  With these identifications and all previous convergences, it is possible to pass to the limit as $m \to \infty$ in each term of the equations of \eqref{problema_P_m_controlado}.


\

\noindent {\bf Energy inequality \eqref{desigualdade_integral_limite}:}

  In order to finish we must prove the energy inequality \eqref{desigualdade_integral_limite}. First we obtain an integral inequality for the solution $(u_m,z_m)$ of \eqref{problema_P_u_m_z_m_controlado}, where we remind that $z_m = \sqrt{v_m + \alpha^2}$, for small enough but fixed $\alpha > 0$, being $(u_m,v_m)$ the solution of \eqref{problema_P_m_controlado}. According to Lemmas \eqref{lemma_estimativa_u_v_m_1_s=1} and \eqref{lemma_estimativa_u_v_m_1_s_intermediario}, $(u_m,z_m)$ satisfies
  \begin{equation} \label{desigualdade_integral_m}
    \begin{array}{l}
      E_m(u_m,z_m)(t_2) + \beta \D{\int_{t_1}^{t_2} \int_{\Omega}}{\norm{\nabla [T^{m}(u_m) + 1]^{s/2}}{}^2 \ dx} \ dt \\[6pt]
      + \dfrac{1}{4} \D{\int_{t_1}^{t_2} \int_{\Omega}}{T^m(u_m)^s \norm{\nabla z_m}{}^2 \ dx} \ dt \\[6pt]
      + \beta \Big ( \D{\int_{t_1}^{t_2} \int_{\Omega}}{\norm{D^2 z_m}{}^2 \ dx} \ dt + \D{\int_{t_1}^{t_2} \int_{\Omega}}{\frac{\norm{\nabla z_m}{}^4}{z_m^2} \ dx} \ dt \Big ) \\[6pt]
      \leq E_m(u_m,z_m)(t_1) + C (\mathcal{K}_1^2) \D{\int_{t_1}^{t_2} \norma{\nabla z_m}{L^2(\Omega)}^2} \ dt + \mathcal{K}_1^2 \D{\int_{t_1}^{t_2} \norma{f}{L^2(\Omega)}^2} \ dt,
    \end{array}
  \end{equation}
  where $E_m$ is given by \eqref{energia_E_m}. Next we collect some convergences that can be obtained from the $m$-independent bounds and the weak and strong convergences proved so far and that will be useful to pass to the limit in the inequality \eqref{desigualdade_integral_m}. Recalling that we denote $z = \sqrt{v + \alpha^2}$, we have, in particular,
  \begin{equation} \label{convergencias_passar_limite_desigualdade_integral}
    \begin{array}{c}
      (u_m + 1)^{s/2} \longrightarrow (u + 1)^{s/2} \mbox{ weakly* in } L^{\infty}(0,T;L^2(\Omega)), \\
      (u_m + 1)^{s/2} \longrightarrow (u + 1)^{s/2} \mbox{ weakly in } L^2(0,T;H^1(\Omega)), \\
      \nabla T^m(u_m)^{s/2} \longrightarrow \nabla u^{s/2} \mbox{ weakly in } L^2(Q), \\
      T^m(u_m)^{s/2} \nabla z_m \longrightarrow u^{s/2} \nabla z \mbox{ weakly in } L^2(Q), \\
      \dfrac{\nabla z_m}{\sqrt{z_m}} \longrightarrow \dfrac{\nabla z }{\sqrt{z}} \mbox{ weakly in } L^4(Q), \\
      D^2 z_m \longrightarrow D^2 z \mbox{ weakly in } L^2(Q), \\
      \nabla z_m \longrightarrow \nabla z \mbox{ strongly in } L^2(Q).
    \end{array}
  \end{equation}
  Recalling that we are dealing with the case $s \in [1,2)$, let
  \begin{equation*}
    E(u,z)(t) = \frac{s}{4} \D{\int_{\Omega}}{g(u(t,x)) \, dx} + \frac{1}{2} \int_{\Omega}{\norm{\nabla z(t,x)}{}^2} dx,
  \end{equation*}
  where
  \begin{equation*}
    g(u) = \left \{ 
    \begin{array}{rl}
      (u+1)ln(u+1) - u, & \mbox{if } s = 1,  \\
      \dfrac{u^s}{s(s-1)}, & \mbox{if } s \in (1,2). 
    \end{array}
    \right.
  \end{equation*}
  Then the following convergence will be also necessary.
  \begin{lemma} \label{lema_convergencia_energia_E_m}
    $E_m(u_m,z_m) \longrightarrow E(u,z)$ in $L^1(0,T)$.
  \end{lemma}
  \begin{proof}[\bf Proof]
    From \eqref{convergencias_passar_limite_desigualdade_integral} we have that $\nabla z_m \to \nabla z$ in $L^2(Q)$ which, in particular, leads us to
    \begin{equation*} 
      \int_{\Omega}{\norm{\nabla z_m(t,x)}{}^2}  dx \longrightarrow \int_{\Omega}{\norm{\nabla z(t,x)}{}^2}  dx \quad \mbox{ in } L^1(0,T).
    \end{equation*}
    Then, it remains to prove that
    \begin{equation} \label{convergencia_diferenca}
      \int_{\Omega}{g_m(u_m) \, dx} \longrightarrow \int_{\Omega}{g(u) \, dx} \quad \mbox{ in } L^1(0,T).
    \end{equation}
    We begin by rewriting $g_m(u_m) - g(u)$ as
    \begin{equation} \label{diferenca_g_m_u_m_g_u}
      g_m(u_m) - g(u) = g_m(u_m) - g_m(u) + g_m(u) - g(u).
    \end{equation}
    For the first difference in \eqref{diferenca_g_m_u_m_g_u}, $g_m(u_m) - g_m(u)$, using that $g_m'$ and $g'$ are monotone increasing functions and that $g_m'(r) \leq g'(r)$ for all $r \geq 0$, we have
    \begin{align*}
      \norm{g_m(u_m) - g_m(u)}{} & = \norm{\int_u^{u_m}{g_m'(\theta) \ d \theta}}{} \\
      & \leq \norm{u_m - u}{} (g_m'(u_m) + g_m'(u)) \\
      & \leq \norm{u_m - u}{} (g'(u_m) + g'(u)).
    \end{align*}
    Then, for $s = 1$, we have
    \begin{equation*}
      \norm{g_m(u_m) - g_m(u)}{} \leq C \norm{u_m - u}{} (ln(u_m + 1) + ln(u + 1))
    \end{equation*}
    and, for $s \in (1,2)$, we have
    \begin{equation*}
      \norm{g_m(u_m) - g_m(u)}{} \leq C \norm{u_m - u}{} (\norm{u_m}{}^{s-1} + \norm{u}{}^{s-1}).
    \end{equation*}
    Considering the case $s \in (1,2)$, since from \eqref{limitacao_u_m_L_5s/3} we have $(\norm{u_m}{}^{s-1} + \norm{u}{}^{s-1})$ bounded in $L^{5s/(3s-3)}(Q)$ and, from \eqref{convergencia_u_m_s_intermediario}, we have $\norm{u_m - u}{} \to 0$ in $L^{5s/(2s+3)}(Q)$, we conclude that
    \begin{equation} \label{convergencia_para_primeira_diferenca}
      g_m(u_m) - g_m(u) \longrightarrow 0 \quad \mbox{ in } L^1(Q).
    \end{equation}
    Considering now the case $s = 1$, from \eqref{limitacao_u_m_L_5s/3} we have $ln(u_m + 1) + ln(u + 1)$ bounded in $L^p(Q)$, for all $p \in [1,\infty)$. Then, analogously to the case $s \in (1,2)$, we use \eqref{convergencia_u_m_s_intermediario} and obtain \eqref{convergencia_para_primeira_diferenca} also for $s = 1$. From \eqref{convergencia_para_primeira_diferenca} we conclude, in particular, that
    \begin{equation} \label{convergencia_primeira_diferenca}
      \int_{\Omega}{g_m(u_m) \ dx} - \int_{\Omega}{g_m(u) \ dx} \longrightarrow 0 \quad \mbox{ in } L^1(0,T).
    \end{equation}
    For the second difference in \eqref{diferenca_g_m_u_m_g_u}, $g_m(u) - g(u)$, we use the Dominated Convergence Theorem. In fact, we write
    \begin{equation*}
      g_m(u) - g(u) = \int_0^u{g_m'(\theta) - g'(\theta) \ d \theta}.
    \end{equation*}
    Using this expression one can verify that
    \begin{equation*}
      g_m(u) - g(u) \longrightarrow 0, \ a.e. \ (t,x) \in Q.
    \end{equation*}
    Next we note that $g_m(u) - g(u)$ is bounded by $2 g(u) \in L^1(Q)$. Therefore we conclude, by using the Dominated Convergence Theorem, that $g_m(u) - g(u) \to 0$ in $L^1(Q)$ and, in particular
    \begin{equation} \label{convergencia_segunda_diferenca}
      \int_{\Omega}{g_m(u) \, dx} - \int_{\Omega}{g(u) \, dx} \longrightarrow 0 \quad\mbox{ in } L^1(0,T).
    \end{equation}
    With \eqref{convergencia_primeira_diferenca} and \eqref{convergencia_segunda_diferenca} we obtain \eqref{convergencia_diferenca}, finishing the proof.
  \end{proof}
  
  \begin{lemma}
    For $s \geq 1$ we have $v \in C_w([0,T];H^1(\Omega))$ and $u \in C_w([0,T];L^s(\Omega))$.
  \end{lemma}
  \begin{proof}[\bf Proof]
    For any $s \geq 1$, we have $v \in L^{\infty}(Q) \subset L^{5/3}(Q)$ and $\partial_t v \in L^{5/3}(Q)$, which implies that $v \in C([0,T];L^{5/3}(\Omega))$. Since $v \in L^{\infty}(0,T;H^1(\Omega))$, by Lemma \ref{lemma_continuidade_fraca} one has $v \in C_w([0,T];H^1(\Omega))$. On the other hand, $u \in L^{5s/(3+s)}(0,T;W^{1,5s/(3+s)}(\Omega))$ and $\partial_t u \in L^{5s/(3+s)}(0,T;(W^{1,5s/(4s-3)}(\Omega))')$. Since from \eqref{convergencias_passar_limite_desigualdade_integral}, $u \in L^{\infty}(0,T;L^s(\Omega))$, we conclude, using Lemma \ref{lemma_continuidade_fraca}, that $u \in C_w([0,T];L^s(\Omega))$.
  \end{proof}
  
  Next we pass to the limit in inequality \eqref{desigualdade_integral_m}. Because of Lemma \ref{lema_convergencia_energia_E_m}, we conclude that, up to a subsequence,
  \begin{equation} \label{convergencia_pontual_E_m}
    E_m(u_m,z_m)(t) \longrightarrow E(u,z)(t) \quad a.e. \ t \in [0,T].
  \end{equation}
  Accounting for the properties of $u_m(0) = u^0_m$ and $v_m(0) = v^0$ we conclude that \eqref{convergencia_pontual_E_m} holds, in particular, for $t = 0$.
  
  Therefore, using the convergences \eqref{convergencias_passar_limite_desigualdade_integral}, the weak lower semicontinuity of the norm (Lemma \ref{lemma_semicontinuidade_fraca_inferior}) and the almost everywhere pointwise convergence \eqref{convergencia_pontual_E_m}, we are able to pass to the limit in \eqref{desigualdade_integral_m} and conclude that, for $a.e.$ $t_1,t_2 \in [0,T]$, with $t_2 > t_1$ (including $t_1 = 0$), we have
  \begin{equation*} 
    \begin{array}{l}
      E(u,z)(t_2) + \beta \D{\int_{t_1}^{t_2} \int_{\Omega}}{\norm{\nabla [u + 1]^{s/2}}{}^2  dx} \, dt + \dfrac{1}{4} \D{\int_{t_1}^{t_2} \int_{\Omega}}{u^s \norm{\nabla z}{}^2  dx} \, dt \\[6pt]
      + \beta \Big ( \D{\int_{t_1}^{t_2} \int_{\Omega}}{\norm{D^2 z}{}^2 dx} \, dt + \D{\int_{t_1}^{t_2} \int_{\Omega}}{\frac{\norm{\nabla z}{}^4}{z^2} dx} \, dt \Big ) \\[6pt]
      \leq E(u,z)(t_1) + C (\mathcal{K}_1^2) \D{\int_{t_1}^{t_2} \norma{\nabla z}{L^2(\Omega)}^2} \ dt + \mathcal{K}_1^2 \D{\int_{t_1}^{t_2} \norma{f}{L^2(\Omega)}^2} \ dt.
    \end{array}
  \end{equation*}
  Accounting for the $m$-independent bound for $\nabla z_m$ given in \eqref{first_estimates_z_m} and the strong convergence of $\nabla z_m$ to $\nabla z$ given in \eqref{convergencias_passar_limite_desigualdade_integral}, we have $\norma{\nabla z}{L^2(Q)} \leq \mathcal{K}_1$. And since $\mathcal{K}_1 = \mathcal{K}_1(\| f\|_{L^q(Q)},\|v^0\|_{W^{2-2/q,q}(\Omega)})$, we conclude that there is other constant $\mathcal{K} = \mathcal{K}(\| f\|_{L^q(Q)},\|v^0\|_{W^{2-2/q,q}(\Omega)})$, which is increasing and continuous with respect to $\norma{f}{L^q(Q)}$, such that, for $a.e.$ $t_1,t_2 \in [0,T]$, with $t_2 > t_1$ (including $t_1 = 0$), we have
  \begin{equation} \label{desigualdade_integral_E}
    \begin{array}{l}
      E(u,z)(t_2) + \beta \D{\int_{t_1}^{t_2} \int_{\Omega}}{\norm{\nabla [u + 1]^{s/2}}{}^2  dx} \, dt + \dfrac{1}{4} \D{\int_{t_1}^{t_2} \int_{\Omega}}{u^s \norm{\nabla z}{}^2  dx} \, dt \\[6pt]
      + \beta \Big ( \D{\int_{t_1}^{t_2} \int_{\Omega}}{\norm{D^2 z}{}^2  dx} \, dt + \D{\int_{t_1}^{t_2} \int_{\Omega}}{\frac{\norm{\nabla z}{}^4}{z^2} dx} \, dt \Big ) \\[10pt]
      \leq E(u,z)(t_1) + \mathcal{K}(\| f\|_{L^q(Q)},\|v^0\|_{W^{2-2/q,q}(\Omega)}).
    \end{array}
  \end{equation}
  To finish, we consider the case $s > 1$. For simplicity, let us write \eqref{desigualdade_integral_E} in terms of the energy $E(u,z)(\cdot)$ and the dissipative term $D(u,z)(t_1,t_2)$, for $a.e.$ $t_1,t_2 \in [0,T]$, with $t_2 > t_1$, as
  \begin{equation*}
    E(u,z)(t_2) + D(u,z)(t_1,t_2) \leq E(u,z)(t_1) + \mathcal{K}(\| f\|_{L^q(Q)},\|v^0\|_{W^{2-2/q,q}(\Omega)}).
  \end{equation*}
  Now, let $t_2 \in (t_1, T]$ and let $\{ t_2^n \}_n$ be a sequence such that $t_2^n \to t_2$ as $n \to \infty$ and such that, for all $n$, we have
  \begin{equation} \label{desigualdade_integral_E_resumida}
    E(u,z)(t_2^n) + D(u,z)(t_1,t_2^n) \leq E(u,z)(t_1) + \mathcal{K}(\norma{f}{L^q(Q)},\|v^0\|_{W^{2-2/q,q}(\Omega)}).
  \end{equation}
  If we take the $\liminf$ in both sides of \eqref{desigualdade_integral_E_resumida}, we obtain
  \begin{equation*}
    \liminf_{n \to \infty}{E(u,z)(t_2^n)} + D(u,z)(t_1,t_2) \leq E(u,z)(t_1) + \mathcal{K}(\norma{f}{L^q(Q)},\|v^0\|_{W^{2-2/q,q}(\Omega)}).
  \end{equation*}
  Then, from the definition of $E$ for $s > 1$ and Lemma \ref{lemma_continuidade_fraca}, we have
  \begin{equation*}
    E(u,z)(t_2) \leq \liminf_{n \to \infty}{E(u,z)(t_2^n)}
  \end{equation*}
  and therefore we conclude that, for $s > 1$, the energy inequality \eqref{desigualdade_integral_E} is satisfied for $a.e.$ $t_1 \in [0,T]$, and for all $t_2 \in (t_1,T]$.


  \subsubsection{\bf Case \texorpdfstring{$\boldsymbol{s \geq 2}$}{s geq 2}}
    \label{subsec:estimativas_u_m_s_geq_2}


\

\noindent {\bf Existence of weak solution to \eqref{problema_P_controlado}:}

    The procedure for the case $s \geq 2$ is slightly different. First we note that, integrating the energy inequality \eqref{estimativa_u_v_m_1_s_geq_2} from Lemma \ref{lemma_estimativa_u_v_m_1_s_geq_2} with respect to $t$, we have
        \begin{equation}
          \nabla T^m(u_m)^{s/2} \mbox{ is bounded in } L^2(Q).
          \label{limitacao_nabla_T^m_u_m_s_geq_2}
        \end{equation}
        We also remind that we defined $g_m'(r) = T^m(r)^{s-1}/(s-1)$, for $s \geq 2$. Then we have
        \begin{align*}
          T^m(r)^s = s \int_0^r{(T^m)'(\theta) T^m(\theta)^{s-1}  d\theta} \leq C s \int_0^r{T^m(\theta)^{s-1}  d\theta} = C s(s - 1) g_m(r).
        \end{align*}
        Therefore it also stems from integrating the energy inequality \eqref{estimativa_u_v_m_1_s_geq_2} with respect to $t$ that
        \begin{equation}
          T^m(u_m)^{s/2} \mbox{ is bounded in } L^{\infty}(0,T;L^2(\Omega)).
          \label{limitacao_T^m_u_m_s_geq_2}
        \end{equation}
        From \eqref{limitacao_T^m_u_m_s_geq_2} and \eqref{limitacao_nabla_T^m_u_m_s_geq_2} we can conclude that
        \begin{equation*}
          T^m(u_m)^{s/2} \mbox{ is bounded in } L^{10/3}(Q),
        \end{equation*}
        that is,
        \begin{equation}
          T^m(u_m) \mbox{ is bounded in } L^{5s/3}(Q).
          \label{limitacao_T^m_u_m_5s/3_s_geq_2}
        \end{equation}
        
        For each fixed $m \in \mathbb{N}$, consider the zero measure set $\mathcal{N} \subset (0,T)$ such that
        \begin{equation*}
          u_m(t^{\ast},\cdot), v_m(t^{\ast},\cdot) \in H^1(\Omega), \ \forall t^{\ast} \in (0,T) \setminus \mathcal{N}.
        \end{equation*}
        Then, for each fixed $t^{\ast} \in (0,T) \setminus \mathcal{N}$, let us consider the sets
        \begin{equation*}
          \{ 0 \leq u_m \leq 1 \} = \Big \{ x \in \Omega \ | \ 0 \leq u_m(t^{\ast},x) \leq 1 \Big \}
        \end{equation*}
        and
        \begin{equation*}
          \{ u_m \geq 1 \} = \Big \{ x \in \Omega \ | \ u_m(t^{\ast},x) \geq 1 \Big \}.
        \end{equation*}
        Now note that, since $s \geq 2$, we have
        \begin{align*}
          & \int_{\Omega} T^m(u_m(t^{\ast},x)^2 \norm{\nabla v_m(t^{\ast},x)}{}^2  dx \\
          & \leq \int_{\{ 0 \leq u_m \leq 1 \}}{\norm{\nabla v_m(t^{\ast},x)}{}^2 dx} + \int_{\{ u_m \geq 1 \}}{T^m(u_m(t^{\ast},x))^s \norm{\nabla v_m(t^{\ast},x)}{}^2 dx} \\
          & \leq \int_{\Omega}{\norm{\nabla v_m(t^{\ast},x)}{}^2 dx} + \int_{\Omega}{T^m(u_m(t^{\ast},x))^s \norm{\nabla v_m(t^{\ast},x)}{}^2 dx}.
        \end{align*}
        The last inequality is valid for all $t^{\ast} \in (0,T) \setminus \mathcal{N}$, then if we integrate in the variable $t$ we obtain
        \begin{align*}
          \int_0^\infty{\int_{\Omega}{T^m(u_m(t,x)^2 \norm{\nabla v_m(t,x)}{}^2  dx} \, dt} & \leq \int_0^\infty{\int_{\Omega}{\norm{\nabla v_m(t,x)}{}^2 dx} \,dt} \\
          & + \int_0^\infty{\int_{\Omega}{T^m(u_m(t,x))^s \norm{\nabla v_m(t,x)}{}^2 dx} \, dt}.
        \end{align*}
        Therefore by \eqref{first_estimates_v_m} and \eqref{limitacao_aux_Delta_v_m_s_intermediario} we can conclude that
        \begin{equation} \label{limitacao_termo_misto_u_v_m_s_geq_2}
          T^m(u_m) \nabla v_m \mbox{ is bounded in } L^2(Q).
      \end{equation}
      
      Now we test the $u_m$-equation of problem \eqref{problema_P_m_controlado} by $u_m$. This gives us
      \begin{align*}
        \dfrac{1}{2} \dfrac{d}{dt} \norma{u_m}{L^2(\Omega)}^2 & + \norma{\nabla u_m}{L^2(\Omega)}^2 = \int_{\Omega}{T^m(u_m) \nabla v_m \cdot \nabla u_m \ dx} \\
        & \leq \dfrac{1}{2} \int_{\Omega}{T^m(u_m)^2 \norm{\nabla v_m}{}^2 dx} + \dfrac{1}{2} \norma{\nabla u_m}{L^2(\Omega)}^2,
      \end{align*}
      hence we have
      \begin{equation*}
        \dfrac{d}{dt} \norma{u_m}{L^2(\Omega)}^2 + \norma{\nabla u_m}{L^2(\Omega)}^2 \leq \int_{\Omega}{T^m(u_m)^2 \norm{\nabla v_m}{}^2 dx}.
      \end{equation*}
      Integrating with respect to $t$, we conclude from \eqref{limitacao_termo_misto_u_v_m_s_geq_2} that
      \begin{equation}
        u_m \mbox{ is bounded in } L^{\infty}(0,T;L^2(\Omega))
        \label{limitacao_u_m_s_geq_2}
      \end{equation}
      and
      \begin{equation}
        \nabla u_m \mbox{ is bounded in } L^2(Q).
        \label{limitacao_nabla_u_m_s_geq_2}
      \end{equation}
      
      Then, if we consider the $u_m$-equation of \eqref{problema_P_m_controlado}, by applying \eqref{limitacao_nabla_u_m_s_geq_2} and \eqref{limitacao_termo_misto_u_v_m_s_geq_2} we conclude that
      \begin{equation}
        \partial_t u_m \mbox{ is bounded in } L^2(0,T;(H^1(\Omega))').
        \label{limitacao_u_m_derivada_tempo_s_geq_2}
      \end{equation}
      
      Considering \eqref{limitacao_aux_Delta_v_m_s_intermediario}, \eqref{first_estimates_v_m} and \eqref{limitacao_T^m_u_m_5s/3_s_geq_2} we conclude from the $v_m$-equation of \eqref{problema_P_m_controlado} that
      \begin{equation*} 
        \partial_t v_m \mbox{ is bounded in } L^2(0,T;L^{3/2}(\Omega)).
      \end{equation*}
      
      Now, using \eqref{limitacao_u_m_s_geq_2}, \eqref{limitacao_nabla_u_m_s_geq_2} and \eqref{limitacao_u_m_derivada_tempo_s_geq_2} we can conclude that there is a subsequence of $\{ u_m \}$, still denoted by $\{ u_m \}$, and a limit function $u$ such that
      \begin{equation*} 
        \begin{array}{c}
          u_m \longrightarrow u \mbox{ weakly* in } L^{\infty}(0,T;L^2(\Omega)), \\
          \nabla u_m \longrightarrow \nabla u \mbox{ weakly in } L^2(Q), \\
          \partial_t u_m \longrightarrow u \mbox{ weakly in } L^2 \Big ( 0,\infty;\big (H^1(\Omega) \big )' \Big ).
        \end{array}
      \end{equation*}
      By applying the compactness result Lemma \ref{lema_Simon}, one has
      \begin{equation*} 
        u_m \longrightarrow u \mbox{ strongly in } L^2(Q).
      \end{equation*}
      Using the Dominated Convergence Theorem and \eqref{limitacao_T^m_u_m_5s/3_s_geq_2} we can also prove that
      \begin{equation*}
        T^m(u_m) \longrightarrow u \mbox{ strongly in } L^p(Q), \ \forall p \in (1,5s/3),
      \end{equation*}
      and using Lemma \ref{lema_convergencia_w_elevado_a_s},
      \begin{equation*}
        T^m(u_m)^s \longrightarrow u^s \mbox{ strongly in } L^p(Q), \ \forall p \in (1,5/3).
      \end{equation*}
      
      From the global in time estimate \eqref{limitacao_T^m_u_m_s_geq_2} we can conclude that, up to a subsequence,
        \begin{equation*}
          T^m(u_m) \to u \mbox{ weakly* in } L^{\infty}(0,T;L^s(\Omega)),
       \end{equation*}
       hence, in particular,
       \begin{equation*}
         u \in L^{\infty}(0,T;L^s(\Omega)).
       \end{equation*}
       
       For $s \geq 2$, if we consider the functions $v_m$, we have the same $m$-independent estimates that we had for $s \in [1,2)$. Then we have the same convergences given in \eqref{convergencia_v_m_s_intermediario}.
       
       Following the ideas of Subsection \ref{subsec:estimativas_u_m_s<2}, we can identify the limits of $T^m(u_m) \nabla v_m$ and $T^m(u_m)^s v_m$ with $u \nabla v$ and $u^s v$, respectively.
       

\

\noindent {\bf Energy inequality \eqref{desigualdade_integral_limite}:}
       
  One can reach it by following the reasoning used in Subsection \ref{subsec:estimativas_u_m_s<2} for $s \in [1,2)$.

\

\section{Existence of an Optimal Control}
  \label{section:existence of optimal control}
  
  In the present section we first prove Theorem \ref{teo_existencia_controle_otimo}, in Subsection \ref{subsec: existencia_controle_otimo}, establishing the existence of solution to the minimization problem \eqref{problema_de_minimizacao}. Afterwards,  we prove Theorem \ref{teo_relacao_entre_problemas_de_minimizacao} in Subsection \ref{subsec: relacao_entre_problemas_de_minimizacao}. 

\subsection{Proof of Theorem \ref{teo_existencia_controle_otimo}}
  \label{subsec: existencia_controle_otimo}

  Since the functional $J$ in \eqref{problema_de_minimizacao} is nonnegative,
  \begin{equation*}
    J_{inf} : = \D{\inf_{(u,v,f) \in S_{ad}^M}} J(u,v,f) \geq 0
  \end{equation*}
  is well defined and there is a minimizing sequence $\{ (u_n,v_n,f_n) \} \subset S_{ad}^M$ such that
  \begin{equation*} 
    \lim_{n \to \infty}{J(u_n,v_n,f_n)} = J_{inf}.
  \end{equation*}
  Next we prove that there is $(\overline{u},\overline{v},\overline{f}) \in S_{ad}^M$, that will be defined as the limit of a subsequence of $\{ (u_n,v_n,f_n) \}_n$, such that $J(\overline{u},\overline{v},\overline{f}) = J_{inf}$.
  
 Since $(u_n,v_n,f_n) \in S_{ad}^M$, we have 
  \begin{equation} \label{problema_P_controlado_n}
    \left\{\begin{array}{l}
      \langle \partial_t u_n, \varphi \rangle + \D{\int_\Omega} \nabla u_n\cdot\nabla \varphi \ dx = \D{\int_\Omega} u_n \nabla v_n \cdot \nabla \varphi \ dx \\[6pt]
      \partial_t v_n - \Delta v_n  = - u_n^s v_n + f_n v_n 1_{\Omega_c}, \\[6pt]
      \partial_\eta u_n  |_{\Gamma}  =  \partial_\eta v_n  |_{\Gamma} = 0, \
      u_n(0)  = u^0, \ v_n(0) = v^0,
    \end{array}\right.
  \end{equation}
   for every $\varphi \in L^{5s/(4s-3)}(0,T;W^{1,5s/(4s-3)}(\Omega))$. Denoting $z_n = \sqrt{v_n + \alpha^2}$, we have 
  \begin{equation} \label{desigualdade_integral_sequencia_minimizante}
    \begin{array}{l}
      E(u_n,z_n)(t_2) + \beta \D{\int_{t_1}^{t_2} \int_{\Omega}}{\norm{\nabla [u_n + 1]^{s/2}}{}^2  dx} \, dt \\[6pt]
      + \dfrac{1}{4} \D{\int_{t_1}^{t_2} \int_{\Omega}}{u_n^s \norm{\nabla z_n}{}^2  dx} \,dt + \beta \Big ( \D{\int_{t_1}^{t_2} \int_{\Omega}}{\norm{D^2 z_n}{}^2  dx} \, dt \\[6pt]
      + \D{\int_{t_1}^{t_2} \int_{\Omega}}{\frac{\norm{\nabla z_n}{}^4}{z^2} dx} \, dt \Big ) \leq E(u_n,z_n)(t_1) + \mathcal{K}(M,\norma{v_0}{W^{2-2/q,q}(\Omega)}).
    \end{array}
  \end{equation}
  Since $(u_n,v_n,f_n) \in S^M_{ad}$, we have
  \begin{equation} \label{limitacao_f_n}
    \norma{f_n}{L^q(Q)} \leq M.
  \end{equation}
  Then, comparing $v_n$ with the solution $w_n$ of \eqref{problema_a_comparar}, with $\tilde{f} = f_n$ and $w^0 = v^0$, yields $0 \leq v_n(t,x) \leq w_n(t,x)$ $a.e.$ $(t,x) \in Q$. From  Lemma \ref{lema_problema_a_comparar} and \eqref{limitacao_f_n}, we obtain $\norma{v_n}{L^{\infty}(Q)} \leq \norma{w_n}{L^{\infty}(Q)} \leq \mathcal{K}_1(M)$ and, in particular, we conclude that there is a constant $C(M) > 0$ such that
  \begin{equation} \label{limitacao_z_n_L_infty}
    0 < \alpha \leq z_n(t,x) \leq C(M), \ a.e. \ (t,x) \in Q.
  \end{equation}
  With the energy inequality \eqref{desigualdade_integral_sequencia_minimizante} and the pointwise bound \eqref{limitacao_z_n_L_infty} at hand, we are able to get the same estimates of Subsection \ref{subsec:estimativas_u_m_v_m_e_passagem_ao_limite} and pass to the limit as $n \to \infty$. In fact, from \eqref{desigualdade_integral_sequencia_minimizante}, we conclude that, for $s \geq 1$, we have the following bounds independently of $n$:
  \begin{equation*}
    \begin{array}{c}
      \nabla z_n \mbox{ is bounded in } L^{\infty}(0,T;L^2(\Omega)) \cap L^4(Q), \\
      u_n^{s/2} \nabla z_n \mbox{ and } D^2 z_n \mbox{ are bounded in } L^2(Q).
    \end{array}
  \end{equation*}
  But using the fact that $z_n = \sqrt{v_n + \alpha^2}$ and \eqref{limitacao_z_n_L_infty} we can conclude that
  \begin{equation*}
    \begin{array}{c}
      \nabla v_n \mbox{ is bounded in } L^{\infty}(0,T;L^2(\Omega)) \cap L^4(Q), \\
      u_n^{s/2} \nabla v_n \mbox{ and } \Delta v_n \mbox{ are bounded in } L^2(Q).
    \end{array}
  \end{equation*}
  From \eqref{desigualdade_integral_sequencia_minimizante} (and by testing the $u_n$-equation of \eqref{problema_P_controlado_n} by $\varphi = 1$, in the case $s = 1$) we also have
  \begin{equation*}
    \begin{array}{c}
      \nabla [u_n + 1]^{s/2} \mbox{ is bounded in } L^2(Q), \\
      (u_n + 1)^{s/2} \mbox{ is bounded in } L^{\infty}(0,T;L^2(\Omega)).
    \end{array}
  \end{equation*}
  Afterwards, using some ideas of Subsection \ref{subsec:estimativas_u_m_v_m_e_passagem_ao_limite} we conclude that for $s \in [1,2)$ we have
  \begin{equation*}
    \begin{array}{c}
      u_n \mbox{ is bounded in } L^{5s/(3 + s)}(0,T;W^{1,5s/(3 + s)}(\Omega)), \\
      \partial_t u_n \mbox{ is bounded in } L^{5s/(3 + s)}(0,T;(W^{1,5s/(4s - 3)}(\Omega))'),
    \end{array}
  \end{equation*}
  for $s \geq 2$ we have
  \begin{equation*}
    \begin{array}{c}
      u_n \mbox{ is bounded in } L^2(0,T;H^1(\Omega)), \\
      \partial_t u_n \mbox{ is bounded in } L^2(0,T;(H^1(\Omega))'),
    \end{array}
  \end{equation*}
  and, for $s \geq 1$, we have
  \begin{equation*}
    \begin{array}{c}
      u_n \mbox{ is bounded in } L^{\infty}(0,T;L^s(\Omega)) \cap L^{5s/3}(Q), \\
      v_n \mbox{ is bounded in } L^{\infty}(Q) \cap L^4(0,T;W^{1,4}(\Omega)) \cap L^2(0,T;H^2(\Omega)), \\
      \partial_t v_n \mbox{ is bounded in } L^{5/3}(Q).
    \end{array}
  \end{equation*}
  In view of these $n$-uniform bounds we follow the reasoning of Subsections \ref{subsec:estimativas_u_m_s<2} and \ref{subsec:estimativas_u_m_s_geq_2} and conclude that, up to a subsequence, there is $(\overline{u}, \overline{v},\overline{f})$ such that, if $s \in [1,2)$, we have
  \begin{equation} \label{convergencia_u_n_minimizante_s<2}
    \begin{array}{c}
      u_n \longrightarrow \overline{u} \mbox{ weakly in } L^{5s/3}(Q) \cap L^{5s/(3+s)}(0,T;W^{1,5s/(3+s)}(\Omega)), \\
      u_n \longrightarrow \overline{u} \mbox{ strongly in } L^p(Q), \ \forall p \in [1,5s/3) \\
      \mbox{and } \partial_t u_n \rightarrow \partial_t \overline{u} \mbox{ weakly in } L^{5s/(3 + s)}(0,T;(W^{1,5s/(4s-3)}(\Omega))'),
    \end{array}
  \end{equation}
  for $s \geq 2$ we have
  \begin{equation} \label{convergencia_u_n_minimizante_s>2}
    \begin{array}{c}
      u_n \longrightarrow \overline{u} \mbox{ weakly in } L^{5s/3}(Q) \cap L^2(0,T;H^1(\Omega)), \\
      u_n \longrightarrow \overline{u} \mbox{ strongly in } L^p(Q), \ \forall p \in [1,5s/3) \\
      \mbox{and } \partial_t u_n \rightarrow \partial_t \overline{u} \mbox{ weakly in } L^2(0,T;(H^1(\Omega))'),
    \end{array}
  \end{equation}
  and, for $s \geq 1$,
  \begin{equation} \label{convergencia_v_n_minimizante}
    \begin{array}{c}
      v_n \rightarrow \overline{v} \mbox{ weakly* in } L^{\infty}(Q) \cap L^{\infty}(0,T;H^1(\Omega)), \\
      v_n \rightarrow \overline{v} \mbox{ weakly in } L^4(0,T;W^{1,4}(\Omega)) \cap L^2(0,T;H^2(\Omega)), \\
      \mbox{and } \partial_t v_n \rightarrow \partial_t \overline{v} \mbox{ weakly in } L^{5/3}(Q)
    \end{array}
  \end{equation}
  and
  \begin{equation} \label{convergencia_f_n_minimizante}
    f_n \rightarrow \overline{f} \mbox{ weakly in } L^q(Q).
  \end{equation}
  With these convergences we can pass to the limit as $n \to \infty$ in \eqref{problema_P_controlado_n} and conclude that $(\overline{u},\overline{v})$ is a weak solution of \eqref{problema_P_controlado} with control $\overline{f}$.
  
  Now we are going to prove that $(\overline{u},\overline{v})$ satisfies the energy inequality \eqref{desigualdade_integral_sequencia_minimizante} and therefore $(\overline{u},\overline{v},\overline{f}) \in S_{ad}^M$. Let $\overline{z} = \sqrt{\overline{v} + \alpha^2}$, if we continue following the ideas of Subsections \ref{subsec:estimativas_u_m_s<2} and \eqref{subsec:estimativas_u_m_s_geq_2} we conclude the convergences 
  \begin{equation} \label{convergencias_passar_limite_desigualdade_integral_sequencia_minimizante}
    \begin{array}{c}
      (u_n + 1)^{s/2} \longrightarrow (\overline{u} + 1)^{s/2} \mbox{ weakly* in } L^{\infty}(0,T;L^2(\Omega)), \\
      (u_n + 1)^{s/2} \longrightarrow (\overline{u} + 1)^{s/2} \mbox{ weakly in } L^2(0,T;H^1(\Omega)), \\
      u_n^{s/2} \nabla z_n \longrightarrow \overline{u}^{s/2} \nabla \overline{z} \mbox{ weakly in } L^2(Q), \\[6pt]
      \dfrac{\nabla z_n}{\sqrt{z_n}} \longrightarrow \dfrac{\nabla \overline{z}}{\sqrt{\overline{z}}} \mbox{ weakly in } L^4(Q), \\[12pt]
      D^2 z_n \longrightarrow D^2 \overline{z} \mbox{ weakly in } L^2(Q), \\
      \nabla z_n \longrightarrow \nabla \overline{z} \mbox{ strongly in } L^2(Q),
    \end{array}
  \end{equation}
  \begin{equation} \label{convergencia_energia_sequencia_minimizante}
    E(u_n,z_n)(t) \rightarrow E(\overline{u},\overline{z})(t) \ a.e. \ t \in (0,T),
  \end{equation}
  and the weak continuity regularity 
  \begin{equation} \label{continuidade_fraca_sequencia_minimizante}
    \overline{v} \in C_w([0,T];H^1(\Omega)) \mbox{ and } \overline{u} \in C_w([0,T];L^s(\Omega)), \mbox{ for } s \geq 1.
  \end{equation}
  Therefore, using \eqref{convergencias_passar_limite_desigualdade_integral_sequencia_minimizante}, \eqref{convergencia_energia_sequencia_minimizante} and \eqref{continuidade_fraca_sequencia_minimizante} we are able to pass to the limit in the energy inequality \eqref{desigualdade_integral_sequencia_minimizante} and conclude that $(\overline{u},\overline{v},\overline{f}) \in S_{ad}^M$.
  
  Finally, we prove that the infimum is attained in $(\overline{u},\overline{v},\overline{f})$. Since $(\overline{u},\overline{v},\overline{f}) \in S_{ad}^M$, we have $J_{inf} \leq J(\overline{u},\overline{v},\overline{f})$. On the other hand, considering again Lemma \ref{lemma_semicontinuidade_fraca_inferior}, the functional $J$ is weakly lower semicontinuous and then
  \begin{equation*}
    J(\overline{u},\overline{v},\overline{f}) \leq \liminf_{n \to \infty}{J(u^n,v^n,f^n)} = J_{inf}.
  \end{equation*}
  Therefore we conclude that there is at least one $(\overline{u},\overline{v},\overline{f}) \in S_{ad}^M$ such that $J(\overline{u},\overline{v},\overline{f}) = J_{inf}$, as we wanted to prove.


\subsection{Proof of Theorem \ref{teo_relacao_entre_problemas_de_minimizacao}}
  \label{subsec: relacao_entre_problemas_de_minimizacao}
  
  Since the functional $J$ in \eqref{problema_de_minimizacao_S_ad^E} is nonnegative,
  \begin{equation*}
    J_{inf}^E : = \D{\inf_{(u,v,f) \in S_{ad}^E}} J(u,v,f) \geq 0
  \end{equation*}
  is well defined and there is a sequence $\{ (u_n,v_n,f_n) \} \subset S_{ad}^E$ such that
  \begin{equation} \label{lim_J(u_n,v_n,f_n)_problema_sem_M}
    \lim_{n \to \infty}{J(u_n,v_n,f_n)} = J_{inf}^E.
  \end{equation}
  Since $(u_n,v_n,f_n) \in S_{ad}^E$, in particular it satisfies the system \eqref{problema_P_controlado_n} and the energy inequality
  \begin{equation} \label{desigualdade_integral_sequencia_minimizante_S_ad}
    \begin{array}{l}
      E(u_n,z_n)(t_2) + \beta \D{\int_{t_1}^{t_2} \int_{\Omega}}{\norm{\nabla [u_n + 1]^{s/2}}{}^2  dx} \ dt + \dfrac{1}{4} \D{\int_{t_1}^{t_2} \int_{\Omega}}{u_n^s \norm{\nabla z_n}{}^2  dx} \ dt \\[6pt]
      + \beta \Big ( \D{\int_{t_1}^{t_2} \int_{\Omega}}{\norm{D^2 z_n}{}^2 \ dx} \ dt + \D{\int_{t_1}^{t_2} \int_{\Omega}}{\frac{\norm{\nabla z_n}{}^4}{z^2} \ dx} \ dt \Big ) \\[10pt]
      \leq E(u_n,z_n)(t_1) + \mathcal{K}(\norma{f_n}{L^q(Q)},\norma{v_0}{W^{2-2/q,q}(\Omega)}).
    \end{array}
  \end{equation}
  
  Following the proof of Theorem \ref{teo_existencia_controle_otimo}, in Subsection \ref{subsec: existencia_controle_otimo}, but this time using \eqref{desigualdade_integral_sequencia_minimizante_S_ad}, we conclude that there is a continuous and increasing function of $\norma{f_n}{L^q(Q)}$, let us denote it by $C(\norma{f_n}{L^q(Q)}) > 0$, such that
  \begin{equation*} 
    \begin{array}{c}
      \norma{u_n}{L^{5s/3}(Q)}, \norma{\nabla u_n}{L^{5s/(3+s)}(Q)} \leq C(\norma{f_n}{L^q(Q)}), \\[6pt]
      \norma{\partial_t u_n}{L^{5s/(3 + s)}(0,T;W^{1,5s/(4s-3)}(\Omega)')} \leq C(\norma{f_n}{L^q(Q)}),
    \end{array}
  \end{equation*}
  for $s \in [1,2)$,
  \begin{equation*} 
    \norma{u_n}{L^{5s/3}(Q)}, \norma{\nabla u_n}{L^2(Q)}, \norma{\partial_t u_n}{L^2(0,T;H^1(\Omega)')} \leq C(\norma{f_n}{L^q(Q)}),
  \end{equation*}
  for $s \geq 2$ and
  \begin{equation*} 
    \begin{array}{c}
      \norma{v_n}{L^{\infty}(Q)}, \norma{\nabla v_n}{L^4(Q)} \leq C(\norma{f_n}{L^q(Q)}) \\[6pt]
      \norma{\nabla v_n}{L^\infty(0,T;L^2(\Omega))\cap L^2(0,T;H^1(\Omega))}, \norma{\partial_t v_n}{L^{5/3}(Q)} \leq C(\norma{f_n}{L^q(Q)})
    \end{array}
  \end{equation*}
  for $s \geq 1$. From the definition of the functional $J$ and \eqref{lim_J(u_n,v_n,f_n)_problema_sem_M} we also conclude that
  \begin{equation*}
    f_n \mbox{ is bounded in } L^q(Q).
  \end{equation*}
  Analogously to Subsection \ref{subsec: existencia_controle_otimo}, from the latter and \eqref{desigualdade_integral_sequencia_minimizante_S_ad} we prove that there is $(\overline{u},\overline{v},\overline{f})\in X_u\times X_v\times L^q(Q)$ such that, up to a subsequence, we have the convergences \eqref{convergencia_u_n_minimizante_s<2}, \eqref{convergencia_u_n_minimizante_s>2}, \eqref{convergencia_v_n_minimizante} and \eqref{convergencia_f_n_minimizante}. These convergences allow us to conclude that $(\overline{u},\overline{v})$ is a weak solution of \eqref{problema_P_controlado} with control $\overline{f}$. Because of the weakly lower semicontinuity of the norm (Lemma \ref{lemma_semicontinuidade_fraca_inferior})
  \begin{equation} \label{estimativa_J_barra}
    J(\overline{u},\overline{v},\overline{f}) \leq \liminf J(u_n,v_n,f_n) = J_{inf}^E.
  \end{equation}
  However, we are not able to prove the $(\overline{u},\overline{v},\overline{f}) \in S_{ad}^E$ and then we can not guarantee that $J_{inf}^E = J(\overline{u},\overline{v},\overline{f})$. In fact, following the ideas of Subsection \ref{subsec: existencia_controle_otimo}, we are able to take the $\liminf$ in \eqref{desigualdade_integral_sequencia_minimizante_S_ad} and, using  that the map $r\in\mathbb{R}_+\mapsto \mathcal{K}(r, \norma{v_0}{W^{2-2/q,q}(\Omega)}) $ is continuous and therefore
  \begin{equation*}
    \liminf_{n \to \infty}{\mathcal{K}(\norma{f_n}{L^q(Q)}, \norma{v_0}{W^{2-2/q,q}(\Omega)})} \leq \mathcal{K}(\liminf_{n \to \infty}{\norma{f_n}{L^q(Q)}}, \norma{v_0}{W^{2-2/q,q}(\Omega)}),
  \end{equation*}
  we obtain
  \begin{equation} \label{desigualdade_integral_sequencia_minimizante_S_ad_limite}
    \begin{array}{l}
      E(\overline{u},\overline{z})(t_2) + \beta \D{\int_{t_1}^{t_2} \int_{\Omega}}{\norm{\nabla [\overline{u} + 1]^{s/2}}{}^2 \ dx} \ dt + \dfrac{1}{4} \D{\int_{t_1}^{t_2} \int_{\Omega}}{\overline{u}^s \norm{\nabla \overline{z}}{}^2 \ dx} \ dt \\[6pt]
      + \beta \Big ( \D{\int_{t_1}^{t_2} \int_{\Omega}}{\norm{D^2 \overline{z}}{}^2 \ dx} \ dt + \D{\int_{t_1}^{t_2} \int_{\Omega}}{\frac{\norm{\nabla \overline{z}}{}^4}{z^2} \ dx} \ dt \Big ) \\[8pt]
      \leq E(\overline{u},\overline{z})(t_1) + \mathcal{K}(\D{\liminf_{n \to \infty}} \norma{f_n}{L^q(Q)},\norma{v_0}{W^{2-2/q,q}(\Omega)}).
    \end{array}
  \end{equation}
  Since by Lemma \ref{lemma_semicontinuidade_fraca_inferior} we have $\norma{\overline{f}}{L^q(Q)} \leq \D{\liminf_{n \to \infty}} \norma{f_n}{L^q(Q)}$, it is not clear that $(\overline{u},\overline{z})$ satisfies \eqref{desigualdade_integral_limite}, that is, we can not guarantee that $(\overline{u},\overline{v},\overline{f}) \in S_{ad}^E$.
  
  On the other hand, we can prove that for $M \geq \frac{q}{\gamma_f} J_{inf}^E$ we have $(\overline{u},\overline{v},\overline{f}) \in S_{ad}^M$. Indeed, because of the convergence \eqref{convergencia_f_n_minimizante}, the weakly lower semicontinuity of the norm (see Lemma \ref{lemma_semicontinuidade_fraca_inferior}), the inequality $\norma{f}{L^q(Q)}^q \leq \frac{q}{\gamma_f} J(u,v,f)$,  
  and \eqref{lim_J(u_n,v_n,f_n)_problema_sem_M}, one has 
  \begin{equation} \label{estimativa_f_barra}
    \norma{\overline{f}}{L^q(Q)}^q \leq \liminf_{n \to \infty} \norma{f_n}{L^q(Q)}^q \leq \frac{q}{\gamma_f} \liminf_{n \to \infty} J(u_n,v_n,f_n) = \frac{q}{\gamma_f} J_{inf}^E.
  \end{equation}
  Then, taking $M^q \geq \frac{q}{\gamma_f} J_{inf}^E$, we have $\norma{\overline{f}}{L^q(Q)} \leq M$. Moreover, from \eqref{desigualdade_integral_sequencia_minimizante_S_ad_limite} and \eqref{estimativa_f_barra} we also have that  $(\overline{u},\overline{z})$ satisfies
  \begin{equation*} 
    \begin{array}{l}
      E(\overline{u},\overline{z})(t_2) + \beta \D{\int_{t_1}^{t_2} \int_{\Omega}}{\norm{\nabla [\overline{u} + 1]^{s/2}}{}^2 dx} \, dt + \dfrac{1}{4} \D{\int_{t_1}^{t_2} \int_{\Omega}}{\overline{u}^s \norm{\nabla \overline{z}}{}^2 dx} \, dt \\[6pt]
      + \beta \Big ( \D{\int_{t_1}^{t_2} \int_{\Omega}}{\norm{D^2 \overline{z}}{}^2 dx} \, dt + \D{\int_{t_1}^{t_2} \int_{\Omega}}{\frac{\norm{\nabla \overline{z}}{}^4}{z^2} dx} \, dt \Big ) \\[10pt]
      \leq E(\overline{u},\overline{z})(t_1) + \mathcal{K}(M,\norma{v_0}{W^{2-2/q,q}(\Omega)}),
    \end{array}
  \end{equation*}
  which allows us to conclude that $(\overline{u},\overline{v},\overline{f}) \in S_{ad}^M$. Therefore, using \eqref{estimativa_J_barra}, 
  \begin{equation*}
    \min_{(u,v,f) \in S_{ad}^M} J(u,v,f) \leq J(\overline{u},\overline{v},\overline{f}) \le \inf_{(u,v,f) \in S_{ad}^E} J(u,v,f),
  \end{equation*}
  as we wanted to prove.



\

  \section*{Statements and Declarations}
  
  \subsection*{Funding}
    This work has been partially supported by Grant PGC2018-098308-B-I00 (MCI/AEI/FEDER, UE). FGG has also been financed in part by the Grant US-1381261 (US/JUNTA/FEDER, UE) and Grant P20-01120 (PAIDI/JUNTA/FEDER,UE).

  \subsection*{Competing interests}

    The authors have no relevant financial or non-financial interests to disclose.

  \subsection*{Author contributions}

    All authors contributed to the study conception and design. All authors read and approved the final manuscript.


\bibliographystyle{amsplain}

\end{document}